# Menger curvature and rectifiability

By J. C. Léger

## Introduction

Let us first introduce some basic definitions needed to better understand this introduction and the rest of the paper. For a Borel set $E \subset \mathbb{R}^n$, we call "total Menger curvature of $E$" the nonnegative number $c(E)$ defined by

$$c^2(E) = \iiint_{E^3} c^2(x,y,z) d\mathcal{H}^1(x) d\mathcal{H}^1(y) d\mathcal{H}^1(z)$$

where $\mathcal{H}^1$ is the 1-dimensional Hausdorff measure in $\mathbb{R}^n$, $c(x,y,z)$ is the inverse of the radius of the circumcircle of the triangle $(x,y,z)$, that is, following the terminology of [6], the Menger curvature of the triple $(x,y,z)$.

A Borel set $E \subset \mathbb{R}^n$ is said to be "purely unrectifiable" if for any Lipschitz function $\gamma : \mathbb{R} \to \mathbb{R}^n$, $\mathcal{H}^1(E \cap \gamma(\mathbb{R})) = 0$ whereas it is said to be rectifiable if there exists a countable family of Lipschitz functions $\gamma_i : \mathbb{R} \to \mathbb{R}^n$ such that $\mathcal{H}^1(E \setminus \cup_i \gamma_i(\mathbb{R})) = 0$. It may be seen from this definition that any 1-set $E$ (that is, $E$ Borel and $0 < \mathcal{H}^1(E) < \infty$) can be decomposed into two subsets

$$E = E_{\text{irr}} \cup E_{\text{rect}}$$

where $E_{\text{irr}}$ is purely unrectifiable and $E_{\text{rect}}$ is rectifiable (see [4]). We can now state the main theorem of this paper.

Theorem 0.1. *If $E \subset \mathbb{R}^n$ is a 1-set and $c^2(E) < \infty$ then $E$ is rectifiable.*

Before going on, I would like to mention that this result was previously proved by G. David in a paper which is to remain unpublished. His construction is a kind of variant of P. Jones' Traveling Salesman Theorem (see [5]) and its main drawback is that it is very difficult to extend it to dimensions higher than 1. The construction given here to prove Theorem 0.1 extends naturally to any dimension, the main problem being to find interesting analytic or geometric criteria for it to hold.

Here is a brief account of the origin and main application of Theorem 0.1.

A compact subset $E$ of $\mathbb{C}$ is said to be removable for the bounded analytic functions if the constants are the only bounded analytic functions on $\mathbb{C} \setminus E$.

A well-known conjecture of Vitushkin (see [8]) stated that a compact 1-set in the plane is removable for the bounded analytic functions if and only if it is purely unrectifiable.

In 1996, P. Mattila, M. Melnikov and J. Verdera (see [6]) used the Menger curvature to prove that the conjecture holds under the additional assumption that the 1-set is Ahlfors-regular.

Recall that a closed subset $E$ of $\mathbb{R}^n$ is said to be Ahlfors-regular if there exists a constant $C \geq 1$ such that for any ball $B$ centered on $E$, of diameter less than the diameter of $E$,

$$C^{-1}\text{diam}B \leq \mathcal{H}^1(E \cap B) \leq C\text{diam}B. \tag{0.1}$$

Their final argument is based on a condition very similar to the condition $c^2(E) < \infty$, although stronger, which is sufficient for a set to be rectifiable. At that time, this sufficient rectifiability condition was known to be valid only for Ahlfors-regular sets.

Since then, G. David has proved that the Vitushkin conjecture in its full strength is true (see [1]; see also [2] as an intermediate step). The structure of his proof is almost the same as the one of [6] although the details are much more complicated due to the lack of the uniform estimate (0.1). His final argument is Theorem 0.1.

Let us show why Theorem 0.1 is not void and why it gives a necessary and sufficient condition for a Borel subset of $\mathbb{R}^n$ to be rectifiable. Considering $\mathbb{R}^2 = \mathbb{C}$, we have the very important relation which is the starting point of the work of Mattila, Melnikov and Verdera,

$$c^2(z_1, z_2, z_3) = \sum_{\sigma \in \mathfrak{S}_3} \frac{1}{(z_{\sigma(1)} - z_{\sigma(3)})\overline{(z_{\sigma(2)} - z_{\sigma(3)})}} \tag{0.2}$$

where the sum ranges over the group $\mathfrak{S}_3$ of permutations of three elements. This relation is not hard to check considering that the law of sines gives

$$\begin{aligned} c(x, y, z) &= 4\frac{\text{Area of the triangle}(x, y, z)}{d(x, y)d(x, z)d(y, z)} \\ &= 2\frac{d(x, L_{y,z})}{d(x, y)d(x, z)} \end{aligned}$$

where $L_{y,z}$ is the line through $y$ and $z$ and $d(.,.)$ is the Euclidean distance in $\mathbb{R}^n$.

The relation (0.2) implies that the $L^2$ boundedness of the Cauchy kernel operator associated to an Ahlfors-regular subset $E$ of $\mathbb{C}$ is equivalent to the fact that there exists a constant $C \geq 1$ such that for any ball $B \subset \mathbb{C}$,

$$c^2(E \cap B) \leq C\text{diam}B.$$

This property turns out to be equivalent to the fact that $E$ is contained in a single Ahlfors-regular curve of the plane by a theorem of G. David and

S. Semmes. The same results are valid in $\mathbb{R}^n$ because in that case $c^2$ is related to the vectorial kernel $\frac{x}{|x|^2}$ for which we have the same $L^2$ boundedness properties on Ahlfors regular curves in $\mathbb{R}^n$. Theorem 0.1 is the non scale-invariant version of these results.

Noticing that a finite collection of Lipschitzian images of $[0,1]$ is bounded and contained in an Ahlfors-regular curve, we deduce from Theorem 0.1 and the countable union feature of the definition of rectifiablity the following characterization:

THEOREM 0.2. *If $E \subset \mathbb{R}^n$ is Borel then $E$ is rectifiable if and only if there exists a countable family of Borel subsets $F_n$ such that $\cup_n F_n = E$, $\mathcal{H}^1(F_n) < \infty$ and $c^2(F_n) < \infty$.*

I would like to end this introduction with a possible higher dimensional analogue of Theorem 0.1. Let $d$ be a positive integer and, for a Borel subset $E \subset \mathbb{R}^n$, set $c^{d+1}(E)$ to be

$$\int_{x\in E}\int_{y_0\in E}\cdots\int_{y_d\in E}\left(\frac{d\,(x,<y_0,\ldots,y_d>)}{d(x,y_0)\ldots d(x,y_d)}\right)^{d+1} d\mathcal{H}^d(y_0)\ldots d\mathcal{H}^d(y_d)d\mathcal{H}^d(x)$$

where $\mathcal{H}^d$ is the $d$-dimensional Hausdorff measure on $\mathbb{R}^n$ and $d(x, < y_0,\ldots,y_d >)$ is the distance between $x$ and the $d$-plane going through the $d+1$ points $y_0,\ldots,y_d$. The quantity $c^{d+1}(E)$ equals $c^2(E)$ when $d=1$.

The interested reader may check that the method presented in this paper applies with only slight modifications to prove:

THEOREM 0.3. *If $\mathcal{H}^d(E) < \infty$ and $c^{d+1}(E) < \infty$ then, up to a set of $\mathcal{H}^d$-measure zero, $E$ is contained in a countable collection of Lipschitzian images of $\mathbb{R}^d$ (i.e. $E$ is $d$-rectifiable).*

The main problem of this result is that we completely lost the connection with boundedness problems on singular integrals (Riesz kernels on surfaces for example) which are our central interest.

I would like to thank G. David very much for the many conversations we had about this problem and his constant support, and Helen Joyce for kindly correcting many English language mistakes.

## 1. First reduction

Theorem 0.1 will follow from two propositions. The second and most important one states roughly that if we have some control on a set $F$ and if $c^2(F)$ is very, very small then 99 percent of $F$ has to be contained in the graph of some Lipschitz function. The first proposition only says that if $E$ satisfies the hypothesis of Theorem 0.1 then there is a nontrivial part of it, $F$, which satisfies the requirements of the second proposition.

PROPOSITION 1.1. *Let $E$ be a set satisfying the hypotheses of Theorem* 0.1, *then for all $\eta > 0$, there exists a subset $F$ of $E$ such that*

(i) *$F$ is compact,*

(ii) $c^2(F) \leq \eta \operatorname{diam} F$,

(iii) $\mathcal{H}^1(F) > \frac{\operatorname{diam} F}{40}$,

(iv) *for all $x \in F$, for all $t > 0$, $\mathcal{H}^1(F \cap B(x,t)) \leq 3t$.*

*Proof.* This is a standard uniformisation procedure as described in [4, p. 17]. As $0 < \mathcal{H}^1(E) < \infty$, we know (see [4, Cor. 2.5]) that, for $\mathcal{H}^1$-almost all $x \in E$,

$$\frac{1}{2} \leq \limsup_{t \to 0} \frac{\mathcal{H}^1(E \cap B(x,t))}{2t} \leq 1. \tag{1.1}$$

Set, for a positive integer $m$,

$$E_m = \left\{ x \in E \text{ such that for all } t \in ]0, \frac{1}{m}[, \quad \mathcal{H}^1(E \cap B(x,t)) \leq 3t \right\}.$$

We have that $E_m \subset E_{m+1}$ and from (1.1), $\mathcal{H}^1(E \setminus \cup_m E_m) = 0$. Hence, there exists $m$ such that $\mathcal{H}^1(E_m) \geq \frac{1}{2}\mathcal{H}^1(E)$ and $c^2(E_m) \leq c^2(E) < \infty$.

Set, for $\tau > 0$,

$$\mathcal{I}(\tau) = \int\int\int_{A(\tau)} c^2(x,y,z) d\mathcal{H}^1(x) d\mathcal{H}^1(y) d\mathcal{H}^1(z)$$

where $A(\tau) = \{(x,y,z) \in E_m^3, d(x,y) < \tau \text{ and } d(x,z) < \tau\}$. As $c^2(E_m) < \infty$, $\mathcal{I}(\tau) \to 0$ when $\tau \to 0$ so that we can find $0 < \tau_0 \leq \frac{1}{2}m$ such that

$$\mathcal{I}(\tau_0) < \frac{\eta \mathcal{H}^1(E_m)}{60 \times 8}.$$

Consider now the family of closed balls

$$\mathcal{V} = \left\{ B(x,\tau), \quad x \in E_m, \quad 0 < \tau < \tau_0, \quad \mathcal{H}^1(E_m \cap B(x,\tau) \geq \frac{\tau}{10} \right\}.$$

From (1.1), $\mathcal{V}$ is a Vitali class for $E_m$ and because of Vitali's covering theorem, (see [4, Th 1.10]), there exists a countable subfamily of $\mathcal{V}$ of disjoint balls $B(x_i, \tau_i)$ such that

$$\mathcal{H}^1(E_m \setminus \cup_i B(x_i, \tau_i)) = 0$$

and

$$\mathcal{H}^1(E_m) \leq 2 \sum_i \tau_i + \frac{1}{2}\mathcal{H}^1(E_m).$$

To be complete, in order to get this conclusion, we should remark that by the definition of $\mathcal{V}$,

$$\sum_i \tau_i \leq 10 \sum_i \mathcal{H}^1(B(x_i,\tau_i)\cap E_m) < \infty.$$

Moreover, we have that

$$\sum_i c^2(B(x_i,\tau_i)\cap E_m) \leq \mathcal{I}(\tau_0) \leq \frac{\eta \mathcal{H}^1(E_m)}{8\times 60},$$

so by setting

$$I_b = \left\{ i : \quad c^2(B(x_i,\tau_i)\cap E_m) \geq \frac{\eta\tau_i}{60} \right\},$$

we get

$$\sum_{i\in I_b} c^2(B(x_i,\tau_i)\cap E_m) \geq \frac{\eta \sum_i \tau_i}{60}.$$

Hence, $\displaystyle\sum_{i\in I_b} \tau_i \leq \frac{\mathcal{H}^1(E_m)}{8}$ and $\displaystyle\sum_{i\notin I_b} \tau_i \geq \frac{\mathcal{H}^1(E_m)}{8}$ since $\displaystyle\sum_i \tau_i \geq \frac{\mathcal{H}^1(E_m)}{4}$.

We can find a ball $B(x_i,\tau_i)$ such that

- $\mathcal{H}^1(B(x_i,\tau_i)\cap E_m) \geq \frac{\tau_i}{10}$,
- $c^2(B(x_i,\tau_i)\cap E_m) \leq \eta\frac{\tau_i}{60}$,
- for any ball $B$ centered on $B(x_i,\tau_i)$, $\mathcal{H}^1(B(x_i,\tau_i)\cap E_m\cap B) \leq \frac{3}{2}\mathrm{diam}B$.

We cannot take $F$ to be $B(x_i,\tau_i)\cap E_m$ because there is no reason for it to be compact. To fix this, we just use the interior regularity property of Hausdorff measure (see [4, Th. 1.6]) to find $F$, a compact subset of $B(x_i,\tau_i)\cap E_m$ such that $\mathcal{H}^1(F) \geq \frac{\tau_i}{20}$ and remark that $\tau_i \leq 20\times 3\times \mathrm{diam}F$ to get the conclusions of Proposition 1.1. □

From now on, we will not need Hausdorff measure, the main proposition being in fact a general statement about some measures on $\mathbb{R}^n$. To make this clearer, let us define the "total Menger curvature" of a Borel measure $\mu$ on $\mathbb{R}^n$ to be the nonnegative number $c(\mu)$ such that

$$c^2(\mu) = \int\!\!\int\!\!\int c^2(x,y,z)\,d\mu(x)d\mu(y)d\mu(z).$$

It is clear that the total Menger curvature of a Borel set $E$ is exactly the total Menger curvature of the measure $\mathcal{H}^1$ restricted to $E$. We will spend most of this article proving the following:

PROPOSITION 1.2. *For any $C_0 \geq 10$, there exists a number $\eta > 0$ such that if $\mu$ is any compactly supported Borel measure on $\mathbb{R}^n$ verifying*

- $\mu(B(0,2)) \geq 1$, $\mu(\mathbb{R}^n \backslash B(0,2)) = 0$,
- *for any ball* $B$, $\mu(B) \leq C_0 \mathrm{diam} B$,
- $c^2(\mu) \leq \eta$,

*then there exists a Lipschitz graph* $\Gamma$ *such that*

$$\mu(\Gamma) \geq \frac{99}{100}\mu(\mathbb{R}^n).$$

*Proof of Theorem* 0.1. Taking this proposition for granted, it is not hard to see that if $E$ satisfies the hypotheses of Theorem 0.1 and if $\mathcal{H}^1(E_{\text{irr}}) > 0$ then $E_{\text{irr}}$ satisfies the same hypothesis as well, so that we can find $F \subset E_{\text{irr}}$ using Proposition 1.1 and applying Proposition 1.2 to $40 \times \mathcal{H}^1$ restricted to a rescaled copy of $F$. We are then able to find a Lipschitz graph $\Gamma$ intersecting $E_{\text{irr}}$ in a set of positive measure. To be precise, we should remark that the $c^2$ function of a set scales like a length and is invariant under isometries: this enables us to rescale the set $F$ to a set of diameter 1 contained in the ball $B(0,2)$. We have a contradiction and Theorem 0.1 is proved. □

From now on, $\mu$ will be a measure satisfiying the hypothesis of Proposition 1.2 and we will note its support $F$. Our duty is to find an adequate coordinate system of $\mathbb{R}^n$ and a Lipschitz function $A : \mathbb{R} \to \mathbb{R}^{n-1}$ whose graph will be the one we are looking for.

These will be defined in Section 3 just after a first investigation on how to handle the geometry of the set $F$ with the little information about $F$ we start with.

Before starting the real technicalities, I would like to point out that we assume that every ball appearing in the following construction is closed. This assumption is needed in order to apply Besicovitch's covering lemma. This is not a serious issue for our construction.

## 2. P. Jones' $\beta$ functions

In this section, we define some functions used to measure how well the support of the measure $\mu$ is approximated by straight lines at a given location and a given scale. We will see that these functions are related to the $c^2$ number provided we are looking at points where the measure $\mu$ does not degenerate too much. To quantify this notion of degeneracy, we need the following density functions.

*Definition* 2.1. For a ball $B$ with center $x \in \mathbb{R}^n$ and radius $t > 0$, we set

$$\delta(B) = \delta(x,t) = \frac{\mu(B(x,t))}{t}$$

and, fixing a number $k_0 \geq 1$,

$$\tilde{\delta}(B) = \tilde{\delta}(x,t) = \sup_{y \in B(x,k_0 t)} \delta(y,t).$$

*Definition* 2.2. Let $k > 1$ be some fixed number. For $x \in \mathbb{R}^n$, $t > 0$ and $D$ a line in $\mathbb{R}^n$, we set

$$\begin{aligned}
\beta_1^D(x,t) &= \frac{1}{t}\int_{B(x,kt)} \frac{d\,(y,D)}{t} d\mu(y), \\
\beta_2^D(x,t) &= \left(\frac{1}{t}\int_{B(x,kt)} \left(\frac{d\,(y,D)}{t}\right)^2 d\mu(y)\right)^{\frac{1}{2}}, \\
\beta_1(x,t) &= \inf_D \beta_1^D(x,t), \\
\beta_2(x,t) &= \inf_D \beta_2^D(x,t).
\end{aligned}$$

$\beta_1^D(x,t)$ and $\beta_2^D(x,t)$ are designed to measure the mean distance from the support of $\mu$ to the line $D$ inside the ball $B(x,kt)$. If $\delta(x,t)$ is too low, this interpretation is not valid so that these numbers make sense only if we keep a uniform lower control on the density function $\delta$. (Recall that we supposed that there is an upper control on the function $\delta$, namely $\delta(B) \leq 2C_0$.) For this purpose we introduce a density threshold, that is, a number $\delta > 0$, and analyze what happens in a ball $B$ satisfying $\delta(B) \geq \delta$.

We begin with a lemma depicting the basic geometrical situation in such balls. It will be of constant use throughout the rest of the construction.

LEMMA 2.3. *There exist constants $C_1 \geq 1$ and $C_1' \geq 1$ depending only on $C_0$ and $\delta$ such that given any ball $B$ satisfying $\delta(B) \geq \delta$, there exist two balls $B_1$ and $B_2$ of radius $\frac{\mathrm{diam} B}{2C_1}$ such that*

(i) *their centers are at least $\frac{12 \mathrm{diam} B}{2C_1}$ apart,*

(ii) $\mu(B \cap B_i) \geq \frac{\mathrm{diam} B}{2C_1'}$.

*Proof.* Without loss of generality we may suppose $B = B(0,1)$. Let $C_1$ and $C_1'$ be two constants to be chosen at the end of the construction and suppose that any pair of closed balls of radius $\frac{1}{C_1}$ centered on $F \cap B(0,1)$ satisfies that either their centers are less than $\frac{12}{C_1}$ apart or one of them satisfies $\mu(B \cap B(0,1)) \leq \frac{1}{C_1'}$.

We apply Besicovitch's covering lemma to the covering of $F \cap B(0,1)$ by the balls $B(x, \frac{1}{C_1})$ with $x \in F \cap B(0,1)$ to get $N$ families $\mathcal{B}_m$ of disjoint balls, $N$ depending only on the ambient dimension $n$ , such that the union of these families is still a covering of $F \cap B(0,1)$. Considering volume, we see that each family contains no more than $(2C_1)^n$ balls. We have

$$\delta \leq \sum_{m=1}^{N} \sum_{B \in \mathcal{B}_m} \mu(B \cap B(0,1)).$$

Hence, there is at least one family $\mathcal{B}_m$ such that

$$\sum_{B \in \mathcal{B}_m} \mu(B \cap B(0,1)) \geq \frac{\delta}{N}.$$

We set

$$\mathcal{G} = \left\{ B \in \mathcal{B}_m, \mu(B \cap B(0,1)) \geq \frac{1}{C_1'} \right\}.$$

By the hypothesis, any ball in $\mathcal{G}$ is contained in a single ball of radius $\frac{15}{C_1}$; hence

$$\sum_{B \in \mathcal{G}} \mu(B \cap B(0,1)) \leq \frac{30C_0}{C_1}$$

because of the upper control on $\delta(B)$.

Moreover,

$$\sum_{B \notin \mathcal{G}} \mu(B \cap B(0,1)) \leq \frac{(2C_1)^n}{C_1'},$$

so that

$$\delta \leq N \left( 2^n \frac{C_1^n}{C_1'} + 30C_0 \frac{1}{C_1} \right),$$

which gives the contradiction when $C_1$ and $C_1'$ are well chosen. □

A first consequence of Lemma 2.3 is a Carleson-like estimate on the $\beta_1$ which is a cousin of the estimates used in [3] to do the "corona construction." This is the construction we will follow in this paper with the numerous modifications needed to handle the case of non-Ahlfors-regular sets. It should be noted that if the measure $\mu$ were the $\mathcal{H}^1$-measure on an Ahlfors-regular set $F$ then the corona construction of [3] would give the Lipschitz graph we are looking for directly. Our main problem here is when we do not have any lower control of the mass of a ball and we will have to show that such situations cannot happen too often.

PROPOSITION 2.4. *There exists a constant $C$ depending on $\delta$, $C_0$, $k$, $k_0$ such that*

$$\int \int_0^\infty \beta_1(x,t)^2 \mathbb{1}_{\{\tilde{\delta}(x,t) \geq \delta\}} \frac{d\mu(x)dt}{t} \leq Cc^2(\mu).$$

*Proof.* We need to define first some local version of $c^2$. For a fixed number $k_1 > 1$, we set, for any ball $B = B(x,t)$,

$$c_{k_1}^2(x,t) = \int\int\int_{\mathcal{O}_{k_1}(x,t)} c^2(u,v,w) d\mu(u) d\mu(v) d\mu(w)$$

with

$$\mathcal{O}_{k_1}(x,t) = \left\{(u,v,w) \in (B(x,k_1t))^3, d(u,v) \geq \frac{t}{k_1}, d(u,w) \geq \frac{t}{k_1}, d(v,w) \geq \frac{t}{k_1}\right\}$$

and

$$c_{k_1}^2 = \iiint_{\mathcal{O}_{k_1}} c^2(x,y,z)d\mu(x)d\mu(y)d\mu(z)$$

with

$$\mathcal{O}_{k_1} = \left\{ (x,y,z) \in (\mathbb{R}^n)^3, \begin{array}{c} \frac{1}{k_1}d(x,y) \leq d(x,z) \leq k_1 d(x,y) \\ \text{and} \\ \frac{1}{k_1}d(x,y) \leq d(y,z) \leq k_1 d(x,y) \end{array} \right\}.$$

A straightforward use of Fubini's Theorem gives

$$\int\int_0^\infty c_{k_1}^2(x,t)\frac{d\mu(x)dt}{t^2} \leq C(k_1)c^2_{2k_1^2} \leq C(k_1)c^2(\mu).$$

Moreover, Hölder's inequality (using the fact that $\delta(y,kt) \leq 2C_0$) gives, for any $y \in \mathbb{R}^n$, for any $t > 0$,

$$\beta_1(y,t)^2 \leq C\beta_2(y,t)^2$$

so that in order to prove Proposition 2.4, we only have to prove,

LEMMA 2.5. *For all $k_0 \geq 1$, all $k \geq 2$, all $\delta > 0$, there exists $k_1 \geq 1$ and $C \geq 1$ such that if $x \in \mathbb{R}^n$, $t > 0$ and $\delta(x,t) \geq \delta$, then for any $y \in B(x,k_0t)$,*

$$\beta_2(y,t)^2 \leq C\frac{c_{k_1}^2(x,t)}{t} \leq C\frac{c_{k_1+k_0}^2(y,t)}{t}.$$

We can apply Lemma 2.3 twice to find three balls $B_1$, $B_2$ and $B_3$ enjoying the same properties as the two balls of Lemma 2.3 with perhaps different numbers $C_1$ and $C_1'$. For each ball $B_i$, set

$$Z_i = \left\{u \in F \cap B_i \cap B(x,t), \right.$$
$$\left. \int\int \mathbb{1}_{\{(u,v,w)\in\mathcal{O}_{k_1}(x,t)\}}c^2(u,v,w)d\mu(v)d\mu(w) \leq C'\frac{c_{k_1}^2(x,t)}{t}\right\},$$

where $C'$ is chosen using Chebichev's inequality, depending on $\delta$, such that $\mu(Z_i) \geq \frac{t}{2C_1'}$.

For $z_1 \in Z_1$, we choose $z_2 \in Z_2$ such that

$$\int \mathbb{1}_{\{(z_1,z_2,w)\in\mathcal{O}_{k_1}(x,t)\}}c^2(z_1,z_2,w)d\mu(w) \leq C''\frac{c_{k_1}^2(x,t)}{t^2},$$

where $C''$ depends on $\delta$.

Let $L$ be the line going through $z_1$ and $z_2$.

If $w \in (F \cap B(x,(k+k_0)t)) \setminus (2B_1 \cup 2B_2)$,

$$c^2(z_1, z_2, w) = \left( \frac{2d\,(w, L)}{d\,(w, z_1)\, d\,(w, z_2)} \right)^2$$

and

$$\frac{t}{k_1} \leq \frac{t}{C_1} \leq d\,(z_i, w) \leq (k+k_0)t \leq k_1 t$$

if $k_1$ is sufficiently large. Hence

$$\int_{B(x,(k+k_0)t)\setminus(2B_1\cup 2B_2)} \left( \frac{d\,(w, L)}{t} \right)^2 d\mu(w) \leq C c^2_{k_1}(x, t).$$

It remains to look at what happens in the ball $2B_i$. Chebichev's inequality shows there exists $z_3 \in Z_3$ such that

$$\begin{aligned} \int \mathbb{1}_{\{(z_1,w,z_3)\in \mathcal{O}_{k_1}(x,t)\}} c^2(z_1, w, z_3) d\mu(w) &\leq C'' \frac{c^2_{k_1}(x,t)}{t^2}, \\ \int \mathbb{1}_{\{(w,z_2,z_3)\in \mathcal{O}_{k_1}(x,t)\}} c^2(w, z_2, z_3) d\mu(w) &\leq C'' \frac{c^2_{k_1}(x,t)}{t^2}, \\ \left( \frac{d\,(z_3, L)}{t} \right)^2 &\leq C'' \frac{c^2_{k_1}(x,t)}{t}. \end{aligned}$$

If $L'$ is the line going through $z_1$ and $z_3$ we get, as before,

$$\int_{2B_2} \left( \frac{d\,(w, L')}{t} \right)^2 d\mu(w) \leq C c^2_{k_1}(x, t).$$

Let $w'$ be the projection of $w$ on $L'$ and $w''$ the projection of $w'$ on $L$. We have

$$\begin{aligned} d\,(w, L)^2 &\leq d\,(w, w'')^2 \\ &\leq 2(d\,(w, w')^2 + d\,(w', w'')^2) \\ &\leq 2(d\,(w, L')^2 + d\,(w', L)^2), \end{aligned}$$

and by Thales Theorem, $d\,(w', L) = d\,(z_3, L)\, \frac{d(z_1,w')}{d(z_1,z_3)}$. Hence, as $d\,(z_1, w') \leq (k+k_0)t$ and $d\,(z_1, z_3) \geq \frac{t}{C_1}$,

$$\left( \frac{d\,(w', L)}{t} \right)^2 \leq C \frac{c^2_{k_1}(x,t)}{t},$$

so that

$$\int_{2B_2} \left( \frac{d\,(w, L)}{t} \right)^2 d\mu(w) \leq C c^2_{k_1}(x, t).$$

The same estimate on the ball $2B_1$ gives the lemma. □

We end this section with a lemma which will be of constant use during the construction of the function $A$. It says roughly that at a given point and a given scale where we have a controlled density and a small $\beta_1$, the lines almost realizing this $\beta_1$ are very close to one another.

LEMMA 2.6. *For all $\delta > 0$, there exist $\varepsilon_0 > 0$ and $C > 1$ such that if $x \in F$, $t > 0$, $\varepsilon < \varepsilon_0$ and $\delta(x,t) \geq \delta$, $\delta(y,t) \geq \delta$, $d(x,y) \leq \frac{k}{2}t$, then for every pair of lines $D_1$ and $D_2$ satisfying*

$$\frac{1}{t}\int_{B(x,kt)} \frac{d(z,D_1)}{t} d\mu(z) \leq 10^4\varepsilon \text{ and } \frac{1}{t}\int_{B(y,kt)} \frac{d(z,D_2)}{t} d\mu(z) \leq 10^4\varepsilon,$$

(i) *for all $w \in D_1$, $d(w,D_2) \leq C\varepsilon(t + d(w,x))$ and for all $w \in D_2$, $d(w,D_1) \leq C\varepsilon(t + d(w,x))$,*

(ii) $\text{angle}(D_1,D_2) \leq C\varepsilon$.

*Proof.* By Lemma 2.3, we may find two balls $B_1$ and $B_2$ of radius $\frac{t}{C_1}$ contained in $B(x,kt)$ and $B(y,kt)$ such that $\mu(B_i) \geq \frac{t}{C_1'}$. These balls are at least $\frac{10t}{C_1}$ apart.

Because of the hypothesis and Chebichev's inequality, there exist $z_1 \in B_1$ and $z_2 \in B_2$ such that $d(z_i,D_j) \leq C\varepsilon t$ for $i = 1,2$ and $j = 1,2$. Let $z_{11}$ be the projection of $z_1$ on $D_1$ and $z_{21}$ the projection of $z_2$ on $D_1$. If $w \in D_1$, then $w = \alpha z_{11} + (1-\alpha)z_{21}$; therefore

$$d(w,D_2) \leq |\alpha|(d(z_{11},D_2) + d(z_{21},D_2)) + d(z_{21},D_2)$$

whereas $d(z_{11},D_2) \leq d(z_{11},z_{12}) \leq d(z_{11},z_1) + d(z_1,z_{12}) \leq C\varepsilon t$ and the same is true for for $d(z_{21},D_2)$. Hence

$$d(w,D_2) \leq C(|\alpha| + 1)\varepsilon t.$$

Moreover, $d(w,z_{21}) = |\alpha| d(z_{11},z_{21})$ and $d(z_{11},z_{21}) \geq \frac{t}{C_1}$ and this gives

$$|\alpha| = \frac{d(w,z_{21})}{d(z_{11},z_{21})} \leq \frac{C_1}{t}(d(w,x) + d(x,z_{21})) \leq C(\frac{d(w,x)}{t} + 1).$$

Hence

$$d(w,D_2) \leq C(d(w,x) + t)\varepsilon.$$

The same argument is valid when $w \in D_2$ and this gives the estimate on the angle between $D_1$ and $D_2$. ∎

## 3. Construction of the Lipschitz graph

The construction of the function $A$ will be done by a stopping time argument similar to the one used in [3]. The main difference here is that we are not allowed to use the "dyadic cube" family of partitions which is a central tool of the construction of [3]. Such a family, enjoying so many nice measure and size features, cannot be expected to exist on a set which is not Ahlfors-regular. This will be fixed by consideration of possibly overlapping balls and use of Besicovitch's covering lemma many times.

3.1. *Construction of the stopping time region.* The main parameter of the construction is the density threshold $\delta$ which we already alluded to. It will be fixed to a value depending only on the ambient dimension $n$. To be precise, we take $\delta = 10^{-10}/N$ where $N$ is the overlap constant appearing in Besicovitch's covering lemma. The other parameters of the construction, namely the numbers $k \geq 10$ and $k_0 \geq 10$ appearing in the definitions of $\beta_1$ and $\tilde{\delta}$, a $\beta_1$-threshold $\varepsilon > 0$, a small angle $\alpha > 0$ and the number $\eta > 0$ will have to be tuned during the construction. Roughly speaking, the $k$'s will be chosen depending only on $\delta$ and $\eta \ll \varepsilon^5 \ll \alpha^{25} \ll 1$. We should recall that $\mu$ is a Borel measure satisfying the hypotheses of Proposition 1.2 and that $F$ is the compact support of $\mu$.

Let us choose a point $x_0 \in F$ and then fix a line $D_0$ such that $\beta_1^{D_0}(x_0, 1) \leq \varepsilon$ which will be the domain of the function $A$. (This is possible because of Lemma 2.5.) Consider now

$$S_{\text{total}} = \left\{ (x,t) \in F \times (0,5), \quad \begin{array}{ll} \text{(i)} & \delta(x,t) \geq \frac{1}{2}\delta \\ \text{(ii)} & \beta_1(x,t) < 2\varepsilon \\ \text{(iii)} & \exists D_{x,t} \text{ s.t. } \left\{ \begin{array}{l} \beta_1^{D_{x,t}}(x,t) \leq 2\varepsilon \\ \text{and} \\ \text{angle}(D_{x,t}, D_0) \leq \alpha \end{array} \right. \end{array} \right\}.$$

We have that $F \times [1,5) \subset S_{\text{total}}$ and $S_{\text{total}}$ is not a stopping time region in the sense of [3] because it is not coherent. This means that if a ball $B$ is in $S_{\text{total}}$ we do not know if larger balls with the same center are also in $S_{\text{total}}$. This property will appear to be crucial in the construction. To correct this, we set, for $x \in F$,

$$h(x) = \sup \left\{ t > 0, \exists y \in F, \exists \tau, \frac{t}{3} \geq \tau \geq \frac{t}{4}, x \in B(y, \frac{\tau}{3}) \text{ and } (y,\tau) \notin S_{\text{total}} \right\}$$

and we set

$$S = \{(x,t) \in S_{\text{total}}, t \geq h(x)\}.$$

*Remark* 3.1. If $(x,t) \in S$ and $t' \geq t$ then $(x,t') \in S$.

This feature of David-Semmes' stopping time regions is called coherence and will be used in the following without much warning. We will often consider $S$ as a set of balls and we will say that $B \in S$ if $B = B(x,t)$ and $(x,t) \in S$. The balls $B(x, h(x))$ belong to $S$. They are called the minimal balls of $S$. We are now ready to cut $F$ in four pieces, one of which will appear to be very nice for what we expect to construct and three others where bad events occur.

Our goal will be to prove that these bad pieces carry only a small part of the measure $\mu$.

*Definition* 3.2 (*A partition of* $F$). Let

$$\mathcal{Z} = \{x \in F, h(x) = 0\},$$

$$F_1 = \left\{ x \in F\backslash\mathcal{Z}, \left\{ \begin{array}{l} \exists y \in F, \exists \tau \in [\frac{h(x)}{5}, \frac{h(x)}{2}], x \in B(y, \frac{\tau}{2}) \\ \text{and} \\ \delta(y,\tau) \le \delta \end{array} \right. \right\},$$

$$F_2 = \left\{ x \in F\backslash(\mathcal{Z} \cup F_1), \left\{ \begin{array}{l} \exists y \in F, \exists \tau \in [\frac{h(x)}{5}, \frac{h(x)}{2}], x \in B(y, \frac{\tau}{2}) \\ \text{and} \\ \beta_1(y,\tau) \ge \varepsilon \end{array} \right. \right\},$$

$$F_3 = \left\{ x \in F\backslash(\mathcal{Z} \cup F_1 \cup F_2), \left\{ \begin{array}{l} \exists y \in F, \exists \tau \in [\frac{h(x)}{5}, \frac{h(x)}{2}], x \in B(y, \frac{\tau}{2}) \\ \text{and} \\ \text{angle}(D_{y,\tau}, D_0) \ge \frac{3}{4}\alpha \end{array} \right. \right\}.$$

*Remark* 3.3. If $x \in F_3$ then for $h(x) \le t \le 100h(x)$, $\text{angle}(D_{x,t}, D_0) \ge \frac{1}{2}\alpha$. To see this, apply Lemma 2.6 to get that $\text{angle}(D_{x,h(x)}, D_{x,t}) \le C\varepsilon$ for each of these $t$ and remember that $\varepsilon \ll \alpha$.

LEMMA 3.4.

$$F = \mathcal{Z} \cup F_1 \cup F_2 \cup F_3$$

*and this union is disjointed.*

*Proof.* Suppose $x \in F\backslash\mathcal{Z}$ so that $h(x) > 0$; then there are sequences $t_n$, $0 < t_n < h(x)$, $t_n \to h(x)$, $y_n \in F$ and $\tau_n$, $\frac{t_n}{4} \le \tau_n \le \frac{t_n}{3}$, such that $x \in B(y_n, \frac{\tau_n}{3})$ and $y_n \notin S_{\text{total}}$ which means that

(1) either $\delta(y_n,\tau_n) < \frac{1}{2}\delta$,

(2) or $\delta(y_n,\tau_n) \ge \frac{1}{2}\delta$ and $\beta_1(y_n,\tau_n) \ge 2\varepsilon$,

(3) or $\delta(y_n,\tau_n) \ge \frac{1}{2}\delta$, $\beta_1(y_n,\tau_n) < 2\varepsilon$ and for any line $\Delta$ such that $\beta_1^{\Delta}(y_n,\tau_n) \le 2\varepsilon$, we have $\text{angle}(\Delta, D_0) > \alpha$.

Because $F$ is compact, we may suppose that $y_n \to y \in F$, $\tau_n \to \tau$ and $(y_n,\tau_n)$ is in case (1) for any $n$ or in case (2) for any $n$ or in case (3) for any $n$. We have that $x \in B(y, \frac{\tau}{3})$ and $\frac{h(x)}{4} \le \tau \le \frac{h(x)}{3}$.

- If $(y_n,\tau_n)$ satisfies (1) for any $n$, then $x \in F_1$.
  Indeed, if $\sigma > 0$ is small, for any sufficiently large $n$, $B(y, \tau - \sigma)$ is contained in the balls $B(y_n,\tau_n)$. Hence

$$\begin{aligned}\mu(B(y,\tau-\sigma)) &\leq \mu(B(y_n,\tau_n))\\ &\leq \frac{1}{2}\delta\tau_n\\ &\leq \delta(\tau-\sigma),\end{aligned}$$

which gives that $\delta(y,\tau-\sigma)\leq\delta$ so that $x\in F_1$.

- If $(y_n,\tau_n)$ satisfies (2) for any $n$ and $x\notin F_1$ then $x\in F_2$.
  Indeed, let $\Delta$ be a line and let $\sigma>0$ be such that $\sigma+\tau\leq\frac{h(x)}{2}$. The ball $B(y,\tau+\sigma)$ contains $B(y_n,\tau_n)$ for any sufficiently large $n$. Hence

$$\begin{aligned}\beta_1^{\Delta}(y,\tau+\sigma) &\geq \left(\frac{\tau_n}{\tau+\sigma}\right)^2\beta_1^{\Delta}(y_n,\tau_n)\\ &\geq 2\left(\frac{\tau_n}{\tau+\sigma}\right)^2\varepsilon\\ &\geq \varepsilon,\end{aligned}$$

  which shows that $x\in F_2$.

- If $(y_n,\tau_n)$ satisfies (3) for any $n$ and $x\notin F_1\cup F_2$ then $x\in F_3$.
  Indeed, let $\Delta$ be a line such that $\beta_1^{\Delta}(y,\tau)\leq\varepsilon$; if $\sigma>0$ is sufficiently small,

$$\begin{aligned}\beta_1^{\Delta}(y_n,\tau_n-\sigma) &\leq \left(\frac{\tau_n-\sigma}{\tau}\right)^2\varepsilon\\ &\leq \frac{3}{2}\varepsilon.\end{aligned}$$

  As there is enough $\mu$-measure of $F$ in $B(y_n,\tau_n-\sigma)$, we can conclude that angle$(\Delta_0,D)\geq\frac{3}{4}\alpha$ which implies that $x\in F_3$ provided $\varepsilon\ll\alpha$. □

As indicated before, we are going to construct a Lipschitz function $A:D_0\to D_0^{\perp}$ such that the set $\mathcal{Z}$ is contained in the graph of $A$. Our goal will be to show that $\mu(\mathcal{Z})\geq\frac{99}{100}\mu(F)$ for an adequate set of parameters. It will then be enough to show that $\mu(F_i)\leq 10^{-6}$ for each $i$, which will be done in Propositions 3.19, 3.5 and 5.9. We can handle the case of $F_2$ now.

Proposition 3.5.

$$\mu(F_2)\leq 10^{-6}.$$

*Proof.* We should remember that $\eta\ll\varepsilon^5$ and that $\delta$ has been chosen once and for all. Now, we remark that if $x\in F_2$, then for any $t\in(h(x),2h(x))$,

$$\beta_1(x,t)\geq\frac{\tau}{t}\beta_1(y,\tau)$$

where $(y,\tau)$ appears in the definition of the fact that $x \in F_2$. We have for such $t$'s,

$$\beta_1(x,t) \geq \frac{\varepsilon}{10}.$$

Now, we have, by Proposition 2.4,

$$\begin{aligned} c^2(F) &\geq \frac{1}{C}\int_F \int_0^\infty 1\!\!1_{\{\tilde\delta(x,t)\geq\delta\}}\beta_1(x,t)^2\frac{d\mu(x)dt}{t} \\ &\geq \frac{1}{C}\int_{F_2}\int_{h(x)}^{2h(x)} \beta_1(x,t)^2\frac{d\mu(x)dt}{t} \\ &\geq \frac{1}{C}\int_{F_2}\int_{h(x)}^{2h(x)} \left(\frac{\varepsilon}{10}\right)^2\frac{d\mu(x)dt}{t} \\ &\geq \frac{1}{C}\varepsilon^2\mu(F_2), \end{aligned}$$

which implies

$$\mu(F_2) \leq \frac{C\eta}{\varepsilon^2} \leq 10^{-6}. \qquad \square$$

In order to construct the function $A$, it is natural to introduce $\pi$, the orthogonal projection onto $D_0$ and $\pi^\perp$, the orthogonal projection onto $D_0^\perp$. We do not have any control of the function $h$ whereas our goal is obviously to control the size of the sets where $h$ is positive and that is why we need to work with some smoothed version of $h$ (see the definition of the function $d$ just below). The second thing we need is some way to associate to each point of $D_0$ some "good" point of $F$; this is the meaning of the function $D$ defined below.

*Definition* 3.6 (*The functions $d$ and $D$*). For $x \in \mathbb{R}^n$, we set

$$d(x) = \inf_{(X,t)\in S}(d\,(X,x)+t)$$

and for $p \in D_0$, we set

$$D(p) = \inf_{x\in\pi^{-1}(p)} d(x) = \inf_{(X,t)\in S}(d\,(\pi(X),p)+t).$$

The following two remarks are easily seen:

*Remark* 3.7. $d$ and $D$ are 1-Lipschitz functions.

*Remark* 3.8. $h(x) \geq d(x)$ and $\mathcal{Z} = \{x \in F, d(x) = 0\}$ because $F$ is closed.

3.2. *Construction of A.* We start with a fundamental lemma which is a first attempt at inverting the projection $\pi : F \to D_0$.

LEMMA 3.9. *There exists a constant $C_2$ such that whenever $x$, $y \in F$ and $t \geq 0$ are such that $d\,(\pi(x),\pi(y)) \leq t$, $d(x) \leq t$, $d(y) \leq t$ then $d\,(x,y) \leq C_2 t$.*

*Proof.* If $t \geq 1$, there is nothing to prove; otherwise, as $d(x) \leq t$ and $d(y) \leq t$, there exist $X$ and $Y$ such that $x \in B(X, 2t) \in S$ and $y \in B(Y, 2t) \in S$.

If $d(x, y) \gg Lt$ where $L$ is some fixed big number, then $d(X, Y) \geq Lt$, $d(\pi(X), \pi(Y)) \leq 5t$, $B_1 = B(X, 2d(X, Y)) \in S$ and $B_2 = B(Y, 2d(X, Y)) \in S$.

Let $D_1$ and $D_2$ be lines associated to $B_1$ and $B_2$ by the definition of $S$. These lines satisfy the hypothesis of Lemma 2.6. Let $B_1' = B(X, 2\varepsilon^{\frac{1}{2}} d(X, Y) + 2t) \in S$ and $B_2' = B(Y, 2\varepsilon^{\frac{1}{2}} d(X, Y) + 2t) \in S$. We have

$$\begin{aligned}
\int_{B_1'} \frac{d(X', D_1)}{d(X, Y)} d\mu(X') &\leq C \frac{d(X, Y)}{\varepsilon^{\frac{1}{2}} d(X, Y)} \frac{1}{2d(X, Y)} \int_{B_1} \frac{d(X', D_1)}{d(X, Y)} d\mu(X') \\
&\leq C\varepsilon^{\frac{1}{2}}
\end{aligned}$$

and

$$\int_{B_2'} \frac{d(Y', D_2)}{d(X, Y)} d\mu(Y') \leq C\varepsilon^{\frac{1}{2}}.$$

Now, by Chebichev, there exist $X' \in B_1'$ and $Y' \in B_2'$ such that $d(X', D_1) \leq C\varepsilon^{\frac{1}{2}} d(X, Y)$ and $d(Y', D_2) \leq C\varepsilon^{\frac{1}{2}} d(X, Y)$.

Considering $X_1'$ the projection of $X'$ on $D_1$, $Y_2'$ the projection of $Y'$ on $D_2$ and $Y_1'$ the projection of $Y_2'$ on $D_1$, we have

$$\begin{aligned}
d(X, X') &\leq (2\varepsilon^{\frac{1}{2}} + \frac{2}{L}) d(X, Y), \\
d(Y, Y') &\leq (2\varepsilon^{\frac{1}{2}} + \frac{2}{L}) d(X, Y), \\
d(X, X_1') &\leq C\varepsilon^{\frac{1}{2}} d(X, Y), \\
d(Y, Y_2') &\leq C\varepsilon^{\frac{1}{2}} d(X, Y), \\
d(Y_1', Y_2') &= d(Y_2', D_1) \leq C\varepsilon d(X, Y)
\end{aligned}$$

by Lemma 2.6. So,

$$\begin{aligned}
d\left(\pi^{\perp}(X), \pi^{\perp}(Y)\right) &\leq d\left(\pi^{\perp}(X), \pi^{\perp}(X')\right) + d\left(\pi^{\perp}(X'), \pi^{\perp}(X_1')\right) \\
&\quad + d\left(\pi^{\perp}(Y_2'), \pi^{\perp}(Y')\right) + d\left(\pi^{\perp}(Y'), \pi^{\perp}(Y)\right) \\
&\quad + d\left(\pi^{\perp}(Y_1'), \pi^{\perp}(Y_2')\right) \\
&\quad + d\left(\pi^{\perp}(X_1'), \pi^{\perp}(Y_1')\right) \\
&\leq C(\varepsilon^{\frac{1}{2}} + \frac{1}{L}) d(X, Y) + 2\alpha d(\pi(X_1'), \pi(Y_1')) \\
&\leq C(\varepsilon^{\frac{1}{2}} + \frac{1}{L}) d(X, Y) + 2\alpha d(\pi(X), \pi(Y)),
\end{aligned}$$

by the same decomposition.

Now if $\varepsilon$ is sufficiently small and $L$ sufficiently large,

$$d\left(\pi^{\perp}(X),\pi^{\perp}(Y)\right)\leq Cd\left(\pi(X),\pi(Y)\right)\leq Ct.$$

Hence $d(X,Y)\leq Ct$, which ends the proof. □

Lemma 3.9 applied with $t=0$ shows that $\pi:\mathcal{Z}\to D_0$ is injective and we can define the function $A$ on $\pi(\mathcal{Z})$ by setting $A(\pi(x))=\pi^{\perp}(x)$ for $x\in\mathcal{Z}$. We should note that a technique similar to the one used in the above proof shows that the function $A:\pi(\mathcal{Z})\to D_0^{\perp}$ is $2\alpha$-Lipschitz, namely

$$d\left(\pi^{\perp}(x),\pi^{\perp}(y)\right)\leq 2\alpha d\left(\pi(x),\pi(y)\right).$$

What remains to do is to extend $A$ to the whole of $D_0$. To that end, we use strictly the same method as in [3] which is a variant of Whitney's extension theorem. Let us choose once and for all a family of dyadic intervals on $D_0$. For $p\in D_0$ such that $p$ is not on the boundary of one of the dyadic intervals and $D(p)>0$, we call $R_p$ the largest dyadic interval containing $p$ and satisfying

$$\operatorname{diam}R_p\leq\frac{1}{20}\inf_{u\in R_p}D(u).$$

The interval $R_p$ does exist because $D(p)>0$. We can now consider the collection of these intervals $R_p$ and relabel it $\{R_i, i\in I\}$. The intervals $R_i$ have disjoint interiors and the family of the $2R_i$'s (here $2R$ is the interval having the same center as $R$ and twice the diameter) is a covering of $D_0\backslash\pi(\mathcal{Z})$. This last fact is due to the fact that $D$ is a 1-Lipschitz function. Using this idea we note:

*Remark* 3.10. If $p\in 10R_i$, $10\operatorname{diam}R_i\leq D(p)\leq 60\operatorname{diam}R_i$.

Indeed, on the one hand, if $u$ is a point in $R_i$, as $D$ is 1-Lipschitz, we have $D(p)\geq D(u)-10\operatorname{diam}R_i\geq 10\operatorname{diam}R_i$; on the other hand, if $u$ is a point of $\tilde{R}_i$, the father of $R_i$ which satisfies $D(u)\leq 20\operatorname{diam}\tilde{R}_i\leq 40\operatorname{diam}R_i$, we have $D(p)\leq D(u)+10\operatorname{diam}\tilde{R}_i\leq 60\operatorname{diam}R_i$.

This gives immediately the following lemma which we will use below.

LEMMA 3.11. (i) *There exists a constant $C$ such that whenever $10R_i\cap 10R_j\neq\emptyset$ then*

$$C^{-1}\operatorname{diam}R_j\leq\operatorname{diam}R_i\leq C\operatorname{diam}R_j.$$

(ii) *For each $i\in I$, there are at most $N$ intervals $R_j$ such that $10R_i\cap 10R_j\neq\emptyset$.*

Notice that we use the same letter $N$ as the one used for Besicovitch's overlap constant: both of them are used in much the same way so that it will not be a problem.

Let us set $U_0 = D \cap B(0,10)$ and $I_0 = \{i \in I, R_i \cap U_0 \neq \emptyset\}$. We claim that there exists a constant $C \geq 1$ such that for any $i \in I_0$, there exists a ball $B_i \in S$ such that

(i) $\mathrm{diam} R_i \leq \mathrm{diam} B_i \leq C \mathrm{diam} R_i$,

(ii) $d(\pi(B_i), R_i) \leq C \mathrm{diam} R_i$.

Indeed, if $p \in R_i$, then there exists $(X_i, t) \in S$ such that $d(p, \pi(X_i)) + t \leq 2D(p) \leq 120 \mathrm{diam} R_i$. Now if $\mathrm{diam} B_i = 2t$ is too small to satisfy (i), we can always go up in $S$ to get $\mathrm{diam} B_i = \mathrm{diam} R_i$ because of Remark 3.1. We let $A_i$ be the affine function $D_0 \to D_0^\perp$ whose graph is $D_i = D_{B_i}$. Also, $A_i$ is Lipschitz of constant $\leq 2\alpha$ (in fact the best constant is less than $\tan \alpha$) because of property (iii) in the definition of $S_{\mathrm{total}}$. We have the following estimates.

Lemma 3.12. *There exists a constant $C$ such that whenever $10R_i \cap 10R_j \neq \emptyset$ then*

(i) $d(B_i, B_j) \leq C \mathrm{diam} R_j$,

(ii) $d(A_i(q), A_j(q)) \leq C\varepsilon \mathrm{diam} R_j$ *for any* $q \in 100 R_j$,

(iii) $|\partial(A_i - A_j)| \leq C\varepsilon$.

*Proof.* For (i), it is enough to apply Lemma 3.9 to the centers of $B_i$ and $B_j$ and to $t = C \mathrm{diam} R_j$. For (ii) and (iii), once we know (i), we can apply Lemma 2.6 provided $k$ is chosen large enough. □

We are now ready to finish the definition of $A$ on $U_0 \backslash \mathcal{Z}$ using a partition of unity. For each $i$, we can find a function $\tilde{\phi}_i \in \mathcal{C}^\infty(D_0)$ such that $0 \leq \tilde{\phi}_i \leq 1$, $\tilde{\phi}_i \equiv 1$ on $2R_i$, $\tilde{\phi}_i \equiv 0$ outside $3R_i$,

$$|\partial \tilde{\phi}_i| \leq \frac{C}{\mathrm{diam} R_i} \text{ and } |\partial^2 \tilde{\phi}_i| \leq \frac{C}{(\mathrm{diam} R_i)^2}.$$

There is then a partition of unity for $V = \bigcup_{i \in I_0} 2R_i$ defined by

$$\phi_i(p) = \frac{\tilde{\phi}_i(p)}{\sum_j \tilde{\phi}_j(p)}.$$

Also,

$$|\partial \phi_i| \leq \frac{C}{\mathrm{diam} R_i} \text{ and } |\partial^2 \phi_i| \leq \frac{C}{(\mathrm{diam} R_i)^2}.$$

Set, for $p \in V$,

$$A(p) = \sum_{i \in I_0} \phi_i(p) A_i(p). \tag{3.1}$$

Since $V \cap \pi(\mathcal{Z}) = \emptyset$ and $U_0 \backslash \pi(\mathcal{Z}) \subset V$, we have just defined $A$ on the whole $U_0$. It remains to show that both definitions glue together to show that $A$ is a $C\alpha$-Lipschitz function.

We show first that $A$ restricted to $2R_j$ is $3\alpha$-Lipschitz. If $p$ and $q$ are two points in $2R_j$,

$$\begin{aligned} d\left(A(p), A(q)\right) &\leq \sum_i \phi_i(p) d\left(A_i(p), A_i(q)\right) \\ &\quad + \sum_i |\phi_i(p) - \phi_i(q)| d\left(A_i(q), A_j(q)\right) \qquad (3.2) \\ &\leq 2\alpha d\left(p,q\right) + C\frac{1}{\mathrm{diam} R_j} d\left(p,q\right) \varepsilon \mathrm{diam} R_j \qquad (3.3) \\ &\leq 3\alpha d\left(p,q\right). \end{aligned}$$

To go from (3.2) to (3.3), note that if $\phi_i(p) - \phi_i(q) \neq 0$ then $3R_i \cap 3R_j \neq \emptyset$ and Lemmas 3.11 and 3.12 apply.

It remains to show that if $p_1 \in \bigcup_{j\in I_0} 2R_j$ and $p_0 \in \pi(\mathcal{Z})$ then

$$d\left(A(p_1), A(p_0)\right) \leq C\alpha d\left(p_1, p_0\right).$$

Let $j \in I_0$ be such that $p_1 \in 2R_j$, pick $p \in R_j$ and let $B_j$ be the corresponding ball, $D_j$ the associated line and $X_j$ a point in $B_j \cap F$ such that $d\left(X_j, D_j\right) \leq C\varepsilon \mathrm{diam} R_j$ (found by Chebichev).

$$\begin{aligned} d\left(A(p_1), A(p_0)\right) &\leq d\left(A(p_1), A(p)\right) \\ &\quad + d\left(A(p), A_j(p)\right) + d\left(A_j(p), A_j(\pi(X_j))\right) \\ &\quad + d\left(A_j(\pi(X_j)), \pi^\perp(X_j)\right) + d\left(\pi^\perp(X_j), A(p_0)\right). \end{aligned}$$

We remark that $D(p_1) \leq d\left(p_1, p_0\right)$ because $D$ is 1-Lipschitz and $D(p_0) = 0$ so that $\mathrm{diam} R_j \leq d\left(p_1, p_0\right)$. Now

$$\begin{aligned} d\left(A(p_1), A(p)\right) &\leq 3\alpha \mathrm{diam} R_j \\ &\leq 3\alpha d\left(p_1, p_0\right), \end{aligned}$$

and

$$\begin{aligned} d\left(A(p), A_j(p)\right) &\leq \sum \phi_i(p) d\left(A_i(p), A_j(p)\right) \\ &\leq C\varepsilon \mathrm{diam} R_j \text{ by Lemma 3.12,} \qquad (3.4) \\ &\leq \alpha d\left(p_1, p_0\right) \text{ because } \varepsilon \ll \alpha, \end{aligned}$$

and

$$\begin{aligned} d\left(A_j(p), A_j(\pi(X_j))\right) &\leq 2\alpha d\left(p, \pi(X_j)\right) \text{ because } A_j \text{ is } 2\alpha\text{-Lipschitz,} \\ &\leq C\alpha \mathrm{diam} R_j \text{ by construction of } B_j, \\ &\leq C\alpha d\left(p_1, p_0\right), \end{aligned}$$

and

$$
\begin{aligned}
d\left(A_j(\pi(X_j)), \pi^{\perp}(X_j)\right) &\leq 2d\left(X_j, D_j\right) \text{ because } \alpha \text{ is small,} \\
&\leq C\varepsilon \mathrm{diam} R_j \text{ by the choice of } X_j, \\
&\leq \alpha d\left(p_1, p_0\right),
\end{aligned}
$$

and, as in the proof of Lemma 3.9 and the proof of the $2\alpha$-Lipschitzness of $A$ on $\pi(\mathcal{Z})$, because $x_0 = p_0 + A(p_0) \in \mathcal{Z}$,

$$
\begin{aligned}
d\left(\pi^{\perp}(X_j), \pi^{\perp}(x_0)\right) &\leq 3\alpha d\left(X_j, x_0\right) \text{ because } \varepsilon^{\frac{1}{2}} \ll \alpha \\
&\leq 10\alpha d\left(\pi(X_j), p_0\right) \text{ because } \alpha \text{ is very small} \\
&\leq 10\alpha(d\left(\pi(X_j), p_1\right) + d\left(p_1, p_0\right)) \\
&\leq C\alpha d\left(p_1, p_0\right).
\end{aligned}
$$

This shows that $A$ is $C\alpha$-Lipschitz.

We end this construction by a last estimate on $A$.

LEMMA 3.13. *There exists a constant $C$ such that if $p \in 2R_j$ then*

$$|\partial^2 A(p)| \leq \frac{C\varepsilon}{\mathrm{diam} R_j}.$$

*Proof.* We have

$$
\begin{aligned}
\partial\partial A &= \partial\partial(\sum_i \phi_i A_i) \\
&= \sum_i (\partial\partial\phi_i) A_i + 2\sum_i \partial\phi_i \partial A_i
\end{aligned}
$$

because $A_i$ is affine. Moreover $\sum_i \partial\phi_i = 0$ so that , for $u \in 2R_j$,

$$|\partial\partial A|(u) \leq \sum_i |\partial\partial\phi_i||A_i - A_j| + 2\sum_i |\partial\phi_i||\partial(A_i - A_j)|.$$

In each of these sums, there are at most $N$ terms; moreover, if $u$ is in the support of $\phi_i$ so that we have $3R_i \cap 3R_j \neq \emptyset$, then $C^{-1}\mathrm{diam} R_i \leq \mathrm{diam} R_j \leq C\mathrm{diam} R_i$. Finally,

$$
\begin{aligned}
|\partial\partial\phi_i| &\leq \frac{C}{(\mathrm{diam} R_i)^2}, \\
|\partial\phi_i| &\leq \frac{C}{\mathrm{diam} R_i}, \\
|A_i - A_j| &\leq C\varepsilon \mathrm{diam} R_i \text{ and} \\
|\partial(A_i - A_j)| &\leq C\varepsilon \text{ by Lemma 3.12}
\end{aligned}
$$

and summing up gives the desired result. □

3.3. *Most of F lies near the graph of A.* The aim of this section is to show that most points of $F$ are at distance less than $\varepsilon^{\frac{1}{2}}d(x)$ from the graph of $A$, which is the thesis of Proposition 3.18.

We set, for $K > 1$,

$$G = \{x \in F \backslash \mathcal{Z}, \forall i, \pi(x) \in 3R_i \Longrightarrow x \notin KB_i\} \cup \{x \in F \backslash \mathcal{Z}, \pi(x) \in \pi(\mathcal{Z})\}.$$

We may remark that if $x \in F \backslash (G \cup \mathcal{Z})$ then there exists $i \in I$ such that $\pi(x) \in 3R_i$ and $x \in KB_i$. Now, if $\pi(x) \in 3R_j$ for $j \neq i$, Lemma 3.11 and construction of the balls $B_i$ guarantee that $\mathrm{diam}R_i$, $\mathrm{diam}R_j$, $\mathrm{diam}B_i$ and $\mathrm{diam}B_j$ are of the same order of magnitude. Lemma 3.12 implies that $x \in K'B_j$ for a $K'$ depending on $K$ and on the other constants appearing in Lemmas 3.11 and 3.12. This shows we could have used “there exists” instead of “for all” in the definition of $G$.

LEMMA 3.14. *If $K$ is large enough, $\mu(G) \leq C\eta$.*

*Proof.* If $x \in G \backslash \pi^{-1}(\mathcal{Z})$ then $\pi(x) \in 3R_i$ for some $i$ and $x \notin KB_i$. Letting $X_i$ be the center of $B_i$, we have

$$\begin{aligned} d(\pi(x), \pi(X_i)) &\leq C\mathrm{diam}R_i \text{ and} \\ d(X_i) &\leq C\mathrm{diam}R_i; \end{aligned}$$

hence, by Lemma 3.9 applied to $x$, $X_i$ and $C\mathrm{diam}R_i$, provided $K$ is large enough, $d(x) \gg \mathrm{diam}B_i$. We can apply the same lemma with $t = d(x)$ to get that $X_i \in B(x, C'd(x))$. Moreover, $d(X_i, x) + \mathrm{diam}B_i \geq d(x)$ so that $d(X_i, x) \geq \frac{d(x)}{2}$ because $\mathrm{diam}B_i \leq \frac{d(x)}{2}$. (See Figure 1.)

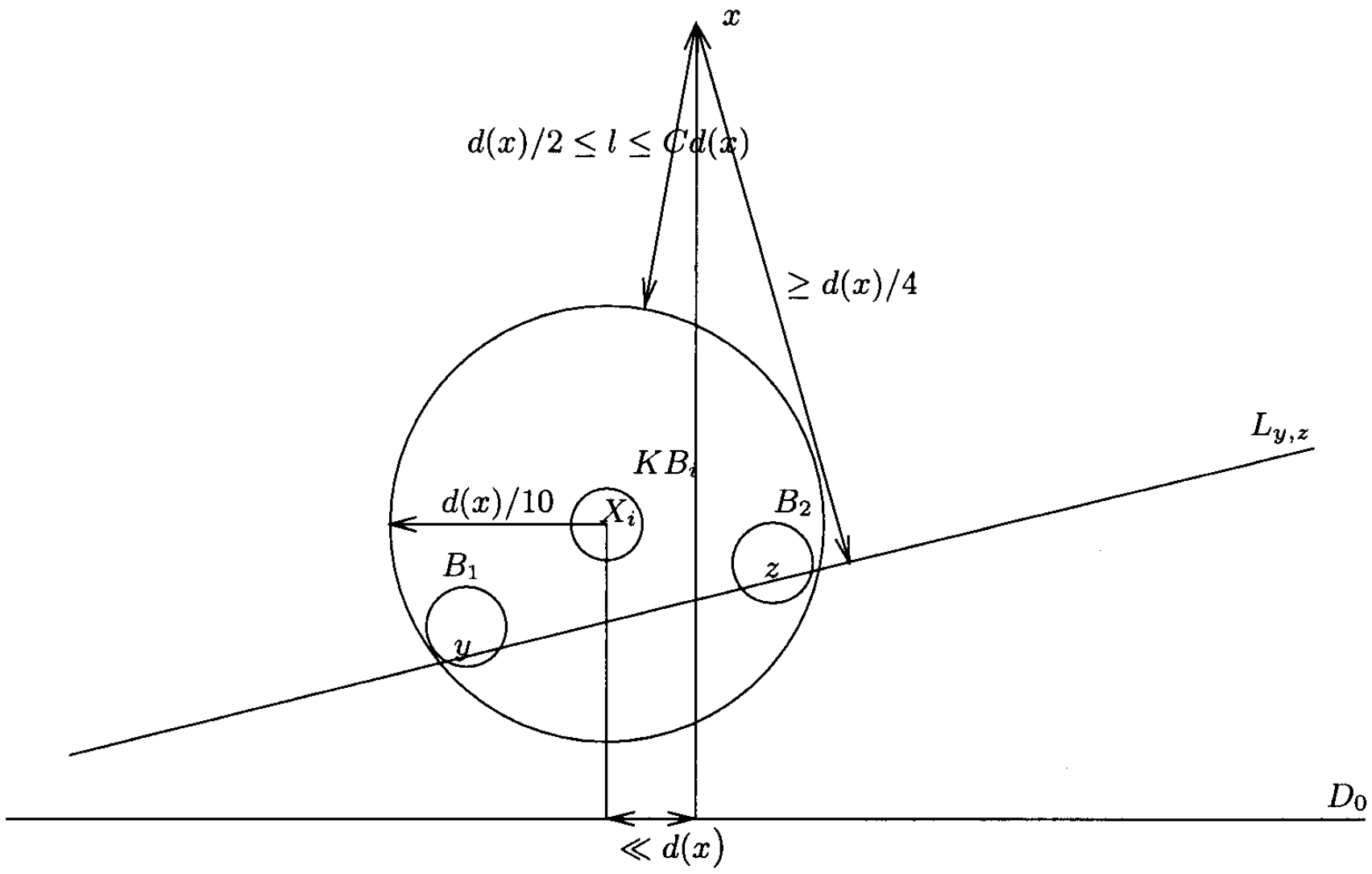

Figure 1. $c^2(x,y,z) \geq \frac{C^{-1}}{d(x)^2}$.

Now, $d(\pi(x), \pi(X_i)) \leq \frac{d(x)}{10}$. The ball $B(X_i, \frac{d(x)}{20})$ belongs to $S$ because it has the same center as $B_i$ and is larger. Using Lemmas 2.3 and 2.6, we find two balls $B_1$ and $B_2$ contained in $B(X_i, \frac{d(x)}{10})$ of radius $\frac{d(x)}{20C_1}$, containing more than $\frac{d(x)}{20C_1'}$ of mass. By Chebichev's inequality, there exist $\tilde{B}_1 \subset B_1$ and $\tilde{B}_2 \subset B_2$ such that $2\mu(\tilde{B}_1) \geq \mu(B_1)$ and $2\mu(\tilde{B}_2) \geq \mu(B_2)$ and for any $y \in \tilde{B}_1$, any $z \in \tilde{B}_2$, the line $L_{y,z}$ going through $y$ and $z$ makes an angle less than $C\varepsilon$ with the line associated with $B(X_i, \frac{d(x)}{10})$. Hence, provided $\varepsilon \ll \alpha$, such a line has an angle less than $2\alpha$ with $D_0$. Hence, if $\alpha$ is small enough, for any $y \in \tilde{B}_1$ and any $z \in \tilde{B}_2$, we have $d(x, L_{y,z}) \geq \frac{d(x)}{4}$ so that

$$\begin{aligned} c^2(x,y,z) &\geq \left(\frac{d(x, L_{y,z})}{d(x,y)\, d(x,z)}\right)^2 \\ &\geq \frac{C^{-1}}{d(x)^2}. \end{aligned}$$

Hence,

$$\begin{aligned} \int\int_{F^2} c^2(x,y,z) d\mu(y) d\mu(z) &\geq \frac{C^{-1}}{d(x)^2} \int_{y \in \tilde{B}_1} \int_{z \in \tilde{B}_2} d\mu(y) d\mu(z) \\ &\geq C^{-1}. \end{aligned}$$

If $x \in G \cap \pi^{-1}(\mathcal{Z})$ we can get the same inequality by reasoning the same way with the point $X = \pi(x) + A(\pi(x)) \in \mathcal{Z}$.

By integrating the inequality over all points $x \in G$, we get

$$\mu(G) \leq C c^2(\mu). \qquad \square$$

Lemma 3.15. *There exists a constant $C_3 \geq 1$ such that for any $x \in F \backslash G$, $C_3^{-1} d(x) \leq D(\pi(x)) \leq d(x)$.*

*Proof.* If $d(x) = 0$, the lemma is obvious; if not, $\pi(x) \in 3R_i$ and $x \in KB_i$ for a given $i$ so that $D(\pi(x)) \geq C^{-1} \mathrm{diam} B_i$; now, $x \in KB_i$ so that $d(x) \leq C \mathrm{diam} B_i$. □

Lemma 3.16. *For any $x \in F$, if $t \geq \frac{d(x)}{10}$,*

$$\int_{B(x,t) \backslash G} d(u, \pi(u) + A(\pi(u)))\, d\mu(u) \leq C\varepsilon t^2.$$

*Proof.* Suppose that $t > 0$, and set

$$I(x,t) = \left\{ i \in I; (2R_i \times D_0^{\perp}) \cap B(x,t) \cap (F \backslash G) \neq \emptyset \right\}.$$

We have

$$\begin{aligned}\int_{B(x,t)\setminus G} & d\left(u,\pi(u)+A(\pi(u))\right)d\mu(u)\\ &\leq \sum_{i\in I(x,t)}\int_{(2R_i\times D_0^\perp)\cap K'B_i} d\left(\pi^\perp(u),A(\pi(u))\right)d\mu(u)\\ &\leq \sum_{i\in I(x,t)}\int_{K'B_i}\sum_j \phi_j(\pi(u))d\left(\pi^\perp(u),A_j(\pi(u))\right)d\mu(u).\end{aligned}$$

Using the facts that

- $d\left(\pi^\perp(u),A_j(\pi(u))\right)\leq 2d\left(u,D_j\right)$ because $\alpha$ is small,
- $\phi_j(\pi(u))\neq 0$ implies, by Lemma 3.11, that $K'B_i\subset kB_j$ provided $k$ is large enough,
- $\beta^{D_j}(B_j)\leq \varepsilon\mathrm{diam}B_j$,
- $\mathrm{diam}B_j$, $\mathrm{diam}R_j$, $\mathrm{diam}R_i$ are of the same order of magnitude,
- there are at most $N$ indices $j$ (see Lemma 3.11 (ii)) such that $\phi_j(\pi(u))\neq 0$,

we get

$$\int_{B(x,t)\setminus G} d\left(u,\pi(u)+A(\pi(u))\right)d\mu(u)\leq C\varepsilon\sum_{i\in I(x,t)}(\mathrm{diam}R_i)^2.$$

Moreover, if $i\in I(x,t)$ then there exists $y\in B(x,t)\cap(F\setminus G)$ such that $\pi(y)\in 2R_i$ so that, because of Remark 3.10 and Lemma 3.15,

$$\mathrm{diam}R_i\leq CD(\pi(y))\leq Cd(y)\leq C(d(x)+t))\leq Ct.$$

Finally, we have

$$\begin{aligned}\int_{B(x,t)\setminus G} d\left(u,\pi(u)+A(\pi(u))\right)d\mu(u) &\leq C\varepsilon t\sum_{i\in I(x,t)}\mathrm{diam}R_i\\ &\leq C\varepsilon t^2\end{aligned}$$

because the cubes $R_i$ are essentially disjoint and are contained in the ball $B(\pi(x),C't)$. □

As we said in the introduction of this section, we want to prove that most points of $F$ are near $\Gamma$, the graph of $A$, which is why we introduce the following definition.

*Definition* 3.17 (A good part of $F$).

$$\tilde{F}=\left\{x\in F\setminus G, d\left(x,\pi(x)+A(\pi(x))\right)\leq \varepsilon^{\frac{1}{2}}d(x)\right\}.$$

We have the following very important proposition.

PROPOSITION 3.18. $\mu(F\backslash\tilde{F}) \leq C\varepsilon^{\frac{1}{2}}$.

*Proof.* We have that $\mu(G) \leq C\eta \leq C\varepsilon^{\frac{1}{2}}$. As

$$F\backslash(\tilde{F}\cup G) \subset \bigcup_{x\in F\backslash(\tilde{F}\cup G)} B(x, \frac{d(x)}{10}),$$

we may use Besicovitch's covering lemma to extract $N$ subfamilies $\mathcal{B}_n$ of disjoint balls from this covering of $F\backslash(\tilde{F}\cup G)$ such that the union of these subfamilies is still a covering of $F\backslash(\tilde{F}\cup G)$. Then

$$\begin{aligned} \varepsilon^{\frac{1}{2}}\mu(F\backslash(\tilde{F}\cup G)) &\leq \int_{F\backslash(\tilde{F}\cup G)} \frac{d\left(u, \pi(u)+A(\pi(u))\right)}{d(u)} d\mu(u) \\ &\leq \sum_{n=0}^{N}\sum_{B\in\mathcal{B}_n}\int_{B\backslash G} \frac{d\left(u, \pi(u)+A(\pi(u))\right)}{d(u)} d\mu(u) \qquad (3.5)\\ &\leq C\sum_{n=0}^{N}\sum_{B\in\mathcal{B}_n}\frac{1}{d(x)}\int_{B\backslash G} d\left(u, \pi(u)+A(\pi(u))\right) d\mu(u) \qquad (3.6)\\ &\leq C\varepsilon\sum_{n=0}^{N}\sum_{B\in\mathcal{B}_n}\mathrm{diam}B \qquad (3.7)\\ &\leq C\varepsilon. \qquad (3.8)\end{aligned}$$

Let us justify these computations.

- To go from (3.5) to (3.6), we note that if $u \in B = B(x, \frac{d(x)}{10})$ then $d(u) \geq \frac{9}{10}d(x)$ because $d$ is 1-Lipschitz.
- To go from (3.6) to (3.7), we apply Lemma 3.16 to the ball $B = B(x, \frac{d(x)}{10})$.
- To go from (3.7) to (3.8), if $B = B(x, \frac{d(x)}{10})$ and $B' = B(y, \frac{d(y)}{10})$ are two balls appearing in the sum, then, provided $C$ is very large (depending on $\delta$),
  - either $\frac{1}{C}B$ and $\frac{1}{C}B'$ have disjoint projections on $D$,
  - or if $2d(x) \geq d(y)$, $B' \subset 4C_2B$ and $\mathrm{diam}B' \geq (2C_3)^{-1}\mathrm{diam}B$.

  Indeed, if the projections of $\frac{1}{C}B$ and $\frac{1}{C}B'$ on $D$ are not disjoint and $d(x) \geq d(y)$, then, applying Lemma 3.9 to $x$, $y$ and $t = d(x)$, we get $d(x, y) \leq C_2 d(x)$ and, by Lemma 3.15, $C_3^{-1}d(x) \leq D(\pi(x)) \leq D(\pi(y)) + \frac{d(x)}{C} \leq d(y) + \frac{d(x)}{C}$ so that $d(y) \geq (2C_3)^{-1}d(x)$.

  Now, we estimate $\sum_{B\in\mathcal{B}_n}\mathrm{diam}B$. Let $B_1$ be a ball in $\mathcal{B}_n$ whose radius is at least half of the maximum radius and let $\mathcal{B}_n^1$ be the family of all the balls $B' \in \mathcal{B}_n$ which satisfy $\pi(\frac{1}{C}B_1) \cap \pi(\frac{1}{C}B') \neq \emptyset$. We can do this operation again with $B_2$, a ball in $\mathcal{B}_n \backslash \mathcal{B}_n^1$ whose radius is at least half of the maximum radius and let $\mathcal{B}_n^2$ be the family of all the balls $B'$ in $\mathcal{B}_n \backslash \mathcal{B}_n^1$ which satisfy $\pi(\frac{1}{C}B_2) \cap \pi(\frac{1}{C}B') \neq \emptyset$. We construct in this way a sequence of balls $B_k$ and a sequence of families $\mathcal{B}_n^k$. Considering volume, we see that each family contains at most $M$ balls where $M$ depends only on $\delta$ because the radii of the balls of each family are of the same order $\sim \mathrm{diam} B_k$ and are contained in a ball which has the same radius. Moreover, we note that

$$\bigcup_k \mathcal{B}_n^k = \mathcal{B}_n$$

because by construction, $\sum_k \mathrm{diam} B_k < \infty$ so that $\mathrm{diam} B_k \to 0$ when $k \to \infty$. Now

$$\begin{aligned}
\sum_{B \in \mathcal{B}_n} \mathrm{diam} B &= \sum_k \sum_{B \in \mathcal{B}_n^k} \mathrm{diam} B \\
&\leq M \sum_k \mathrm{diam} B_k \\
&\leq MC \sum_k \mathrm{diam} \pi(\frac{1}{C} B_k) \\
&\leq C
\end{aligned}$$

because the projections of the balls are disjoint and are contained in $B(0, 10) \cap D_0$. □

We can now estimate $\mu(F_1)$, where, as we recall,

$$F_1 = \left\{ x \in F \backslash \mathcal{Z}, \exists y \in F, \exists \tau \in [\frac{h(x)}{5}, \frac{h(x)}{2}], x \in B(y, \frac{\tau}{2}) \text{ and } \delta(y, \tau) \leq \delta \right\}.$$

Proposition 3.19. $\mu(F_1) \leq 10^{-6}$.

*Proof.* Since

$$F_1 \cap \tilde{F} \subset \bigcup_{x \in F_1 \cap \tilde{F}} B(x, \frac{h(x)}{8}),$$

by Besicovitch's covering lemma, we may extract from this covering $N$ subfamilies $\mathcal{B}_n$ of disjoint balls such that their union is still a covering of $F_1 \cap \tilde{F}$. Notice that, by the construction of $\tilde{F}$, these balls are almost aligned on the graph $\Gamma$. (See Figure 2.) We have then

$$\mu(F_1 \cap \tilde{F}) \leq \sum_{n=0}^{N} \sum_{B \in \mathcal{B}_n} \mu(B) \tag{3.9}$$

$$\leq 10\delta \sum_{n=0}^{N} \sum_{B \in \mathcal{B}_n} \operatorname{diam} B \tag{3.10}$$

$$\leq 1000 N \delta. \tag{3.11}$$

We now justify these computations.

- To go from (3.9) to (3.10), we use the definition of $F_1$,
- To go from (3.10) to (3.11), we note that for a ball $B$ appearing in the sum, we have, provided $\varepsilon$ and $\alpha$ are small enough, $\mathcal{H}^1(\Gamma \cap B) \geq \frac{\operatorname{diam} B}{10}$, so that $\sum_{B \in \mathcal{B}_n} \operatorname{diam} B \leq 10\mathcal{H}^1(\Gamma \cap B(0,2)) \leq 40$.

Having chosen $\delta \leq \frac{10^{-10}}{N}$ (so that $\mu(F_1 \cap \tilde{F}) \leq 10^{-7}$) and $\varepsilon$ very small (so that $\mu(F \backslash \tilde{F}) \leq 10^{-7}$), we obtain the control we sought. □

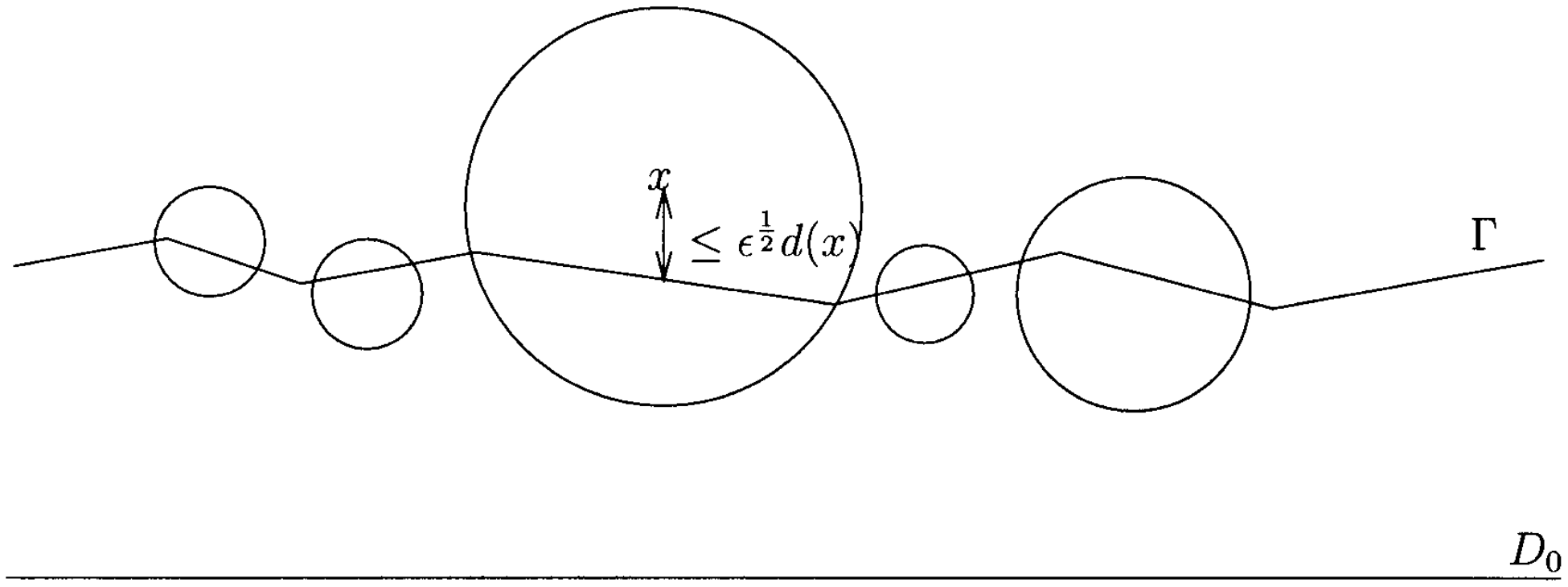

Figure 2. The balls are aligned on $\Gamma$, the graph of $A$.

## 4. The $\gamma$ function of $A$

For $p \in D_0 \cap B(0,10)$ and $t > 0$, set

$$\gamma(p,t) = \inf_a \frac{1}{t} \int_{B(p,t) \cap D_0} \frac{|A(u) - a(u)|}{t} du,$$

where the inf is taken over all affine functions $a : D_0 \to D_0^{\perp}$ and

$$\tilde{\gamma}(p,t) = \inf_M \frac{1}{t} \int_{B(p,t) \cap D_0} \frac{d\,(u + A(u), M)}{t} du,$$

where the infimum is taken over all lines $M$. As the Lipschitz constant of $A$ may be chosen small enough,

$$\frac{1}{2}\tilde{\gamma}(p,t) \le \gamma(p,t) \le 2\tilde{\gamma}(p,t).$$

These $\gamma$ functions measure the approximation of the function $A$ by affine functions and the approximation of the graph of the function $A$ by lines. They are very similar to the $\beta_1$ function and the goal of this part is to get a control on $\gamma$ similar to the one we got on $\beta_1$ in Proposition 2.4, namely

Proposition 4.1.

$$\begin{aligned}\int_0^2 \int_{U_0} \gamma(p,t)^2 \frac{dpdt}{t} &\le C\varepsilon^2 + Cc^2(\mu) \\ &\le C\varepsilon^2\end{aligned}$$

*where $C$ does not depend on $\alpha$.*

We will use this estimate in the next part to show that the function $A$ cannot oscillate too much which would be the case if $F_3$ were too large.

Lemma 4.2.

$$\sum_{i\in I_0} \int_0^{\mathrm{diam} R_i} \int_{R_i} \gamma(p,t)^2 \frac{dpdt}{t} \le C\varepsilon^2.$$

*Proof.* By Taylor's formula, $\gamma(p,t) \le Ct \sup_{u\in B(p,t)} |\partial^2 A(u)|$ so that by Lemma 3.13 (because $u \in 2R_i$),

$$\begin{aligned}\sum_{i\in I_0} \int_0^{\mathrm{diam} R_i} \int_{R_i} \gamma(p,t)^2 \frac{dpdt}{t} &\le C\varepsilon^2 \sum_{i\in I_0} \frac{1}{(\mathrm{diam} R_i)^2} \int_0^{\mathrm{diam} R_i} \int_{R_i} t dpdt \\ &\le C\varepsilon^2 \sum_{i\in I_0} \int_{R_i} dp \\ &\le C\varepsilon^2.\end{aligned}$$

To complete our comparison program, it remains to estimate

$$\int_0^2 \int_{\pi(\mathcal{Z})} \gamma(p,t)^2 \frac{dpdt}{t}$$

and

$$\sum_{i\in I_0} \int_{\mathrm{diam} R_i}^2 \int_{R_i} \gamma(p,t)^2 \frac{dpdt}{t}.$$

Therefore we need an estimate of $\gamma(p,t)$ when $t > \frac{D(p)}{60}$. We fix $p$ and $t$ satisfying this relation. Hence, there exists $(\tilde{X},T) \in S$, $\tilde{X}$ not depending on $t$ such that

- $d\left(\pi(\tilde{X}), p\right) \leq Ct$,
- $T = Ct$.

If $X \in B(\tilde{X}, t) \cap F$, we have $d(X) \leq t + T \leq Ct$. Let $D_{p,t}$ be a line such that

$$\beta_1^{D_{p,t}}(X,t) \leq 2\beta_1(X,t).$$

Now, $I(p,t) = \{i \in I_0, R_i \cap B(p,t) \neq \emptyset\}$. Then,

$$\begin{aligned}
\gamma(p,t) &\leq 2\frac{1}{t}\int_{B(p,t)} \frac{d\left(u + A(u), D_{p,t}\right)}{t} du \\
&\leq 2\frac{1}{t}\int_{B(p,t)\cap\pi(\mathcal{Z})} \frac{d\left(u + A(u), D_{p,t}\right)}{t} du \\
&\quad +2\sum_{i\in I(p,t)} \frac{1}{t}\int_{B(p,t)\cap R_i} \frac{d\left(u + A(u), D_{p,t}\right)}{t} du \\
&= a + \sum_{i\in I(p,t)} a_i.
\end{aligned}$$

We estimate $a$ first. On one end, if $x = u + A(u) \in \mathcal{Z}$ then $d(x) = 0$; on the other end $d\left(\pi(x), \pi(X)\right) \leq Ct$ and $d(X) \leq Ct$ so that $d\left(x, X\right) \leq Ct$ by Lemma 3.9. Hence, as we may push the integral on $\pi(\mathcal{Z})$ on $\mathcal{Z}$ by the parametrization $A$,

$$\begin{aligned}
a &\leq C\frac{1}{t^2}\int_{\mathcal{Z}\cap B(X,Ct)} d\left(x, D_{p,t}\right) d\mathcal{H}^1(x) \\
&\leq C\beta_1(X,t). \qquad (4.1)
\end{aligned}$$

It is worth noticing that to go from the integral against $d\mathcal{H}^1$ to the integral against $d\mu$, we use the fact that for any ball $B$ centered on $\mathcal{Z}$, $2\mu(B) \geq \delta \mathrm{diam} B$. The definition of $\mathcal{H}^1$ and some covering argument implies then that for any function $f$ continuous on $F$,

$$\int_{\mathcal{Z}} f d\mathcal{H}^1 \leq C \int f d\mu.$$

Next, estimating the $a_i$'s, we have

$$\begin{aligned}
a_i &\leq \frac{1}{t}\int_{B(p,t)\cap R_i} \frac{d\left(u + A(u), D_i\right)}{t} du \\
&\quad + \frac{\mathrm{diam} R_i}{t} \sup\left\{\frac{d\left(w, D_{p,t}\right)}{t}, w \in D_i, d\left(w, B_i\right) \leq C \mathrm{diam} R_i\right\} \qquad (4.2)
\end{aligned}$$

because $d\left(u + A(u), D_{p,t}\right) \leq d\left(u + A(u), D_i\right) + d\left(w, D_{p,t}\right)$ where $w$ is the orthogonal projection of $u + A(u)$ on $D_i$.

Moreover, as $D_i$ is the graph of $A_i$, we have, by Lemma 3.12 (see inequality (3.4) as well),

$$\frac{1}{t}\int_{B(p,t)\cap R_i} \frac{d\left(u + A(u), D_i\right)}{t} du \leq C\varepsilon\left(\frac{\mathrm{diam} R_i}{t}\right)^2. \qquad (4.3)$$

LEMMA 4.3.

$$\sup\left\{\frac{d\,(w,D_{p,t})}{t}, w\in D_i, d\,(w,B_i)\le C\mathrm{diam}R_i\right\}$$
$$\le C\varepsilon\frac{\mathrm{diam}R_i}{t}+C\frac{1}{t}\left(\frac{1}{\mathrm{diam}R_i}\int_{2B_i} d\,(z,D_{p,t})^{\frac{1}{3}}\,d\mu(z)\right)^3.$$

*Proof.* Let $B_1$ and $B_2$ be two balls given by Lemma 2.3 applied to the ball $B_i$ and set, for $k=1,2$,

$$Z_k=\left\{z_k\in B_k\cap F, d\,(z_k,D_i)\le C'\varepsilon\mathrm{diam}R_i\right\},$$

where $C'$ is chosen in order to have $\mu(Z_k)\ge\dfrac{\mathrm{diam}R_i}{2000C_1'}$.

If $z_1\in Z_1$ and $z_2\in Z_2$ and if $z_1'$ and $z_2'$ are their projections on $D_i$, we have, provided $\varepsilon$ is very small, that $d\,(z_1',z_2')\ge C^{-1}\mathrm{diam}R_i$. If $w\in D_i$ is such that $d\,(w,B_i)\le C\mathrm{diam}R_i$, we have $w=\sigma z_1'+(1-\sigma)z_2'$ for some $\sigma$ such that $|\sigma|\le d\,(w,z_2')\,d\,(z_1',z_2')^{-1}\le C''$. If $\tilde w$, $\tilde z_1$, $\tilde z_2$ are the projections of $w$, $z_1'$ and $z_2'$ on $D_{p,t}$, we have then

$$\begin{aligned} d\,(w,\tilde w) &\le |\sigma| d\,(z_1',\tilde z_1)+(1+|\sigma|)d\,(z_2',\tilde z_2)\\ &\le C''(d\,(z_1,D_{p,t})+d\,(z_2,D_{p,t}))+C'''\varepsilon\mathrm{diam}R_i. \end{aligned}$$

Hence, after some cube and cubic root manipulations, by integrating on $B_1\cup B_2$, we get

$$\mu((B_1\cup B_2)\cap F)^3\left(d\,(w,D_{p,t})-C'''\varepsilon\mathrm{diam}R_i\right)\le C\left(\int_{B_1\cup B_2} d\,(z,D_{p,t})^{\frac{1}{3}}\,d\mu(z)\right)^3,$$

so that

$$\frac{d\,(w,D_{p,t})}{t}\le C\varepsilon\frac{\mathrm{diam}R_i}{t}+\frac{C}{t}\left(\frac{1}{\mathrm{diam}R_i}\int_{2B_i} d\,(z,D_{p,t})^{\frac{1}{3}}\,d\mu(z)\right)^3. \qquad \square$$

LEMMA 4.4.

$$\begin{aligned} \sum_{i\in I(p,t)}\frac{\mathrm{diam}R_i}{t}\times\frac{1}{t}\left(\frac{1}{\mathrm{diam}R_i}\int_{2B_i} d\,(z,D_{p,t})^{\frac{1}{3}}\,d\mu(z)\right)^3 &\\ \le\ & \frac{C}{t^2}\int_{\bigcup_{i\in I(p,t)}2B_i} d\,(z,D_{p,t})\,d\mu(z)\\ \le\ & C\beta_1(X,t). \end{aligned}$$

*Proof.* For $i\in I(p,t)$, we set

$$\begin{aligned} J(i) &= \left\{j\in I(p,t),\mathrm{diam}B_j\le\mathrm{diam}B_i \text{ and } 2B_i\cap 2B_j\ne\emptyset\right\},\\ N_i(x) &= \sum_{j\in J(i)}1\!\!1_{2B_j}(x). \end{aligned}$$

For $x \in F$ and $k$ an integer, let $B_i$ be a ball of maximal diameter such that $x \in 2B_i$ and $N_i(x) = k$. If $B_j$ is another ball satisfying these properties, save for the maximality, we have $\text{diam}B_i = \text{diam}B_j$ because if this were not the case, we would have $N_j(x) < N_i(x)$. Now, $R_j \subset CR_i$, these dyadic cubes are disjoint and their sizes are comparable; hence, there are at most $C$ of such balls $B_j$, $C$ not depending on $k$. Hence,

$$\begin{aligned}\sum_{i \in I_0} \mathbb{1}_{2B_i}(x) N_i(x)^{-2} &= \sum_k \frac{1}{k^2} \left( \sum_{i \in I_0, N_i(x) = k} \mathbb{1}_{2B_i}(x) \right) \\ &\le C \sum_k \frac{1}{k^2} \\ &\le C\end{aligned}$$

and

$$\begin{aligned}\int_{2B_i \cap F} N_i(x) d\mu(x) &\le \sum_{j \in J(i)} \mu(2B_i \cap F) \\ &\le C \sum_{j \in J(i)} \text{diam}R_j \\ &\le C\text{diam}R_i\end{aligned}$$

because the dyadic cubes $R_j$ are disjoint, of comparable sizes and are within distance $C\text{diam}R_i$ from $R_i$.

By Hölder's inequality, we get

$$\begin{aligned}&\left( \frac{1}{\text{diam}R_i} \int_{2B_i} d\left(z, D_{p,t}\right)^{\frac{1}{3}} N_i(z)^{\frac{-2}{3}} N_i(z)^{\frac{2}{3}} d\mu(z) \right)^3 \\ &\quad\le \left( \frac{1}{\text{diam}R_i} \int_{2B_i} d\left(z, D_{p,t}\right) N_i(z)^{-2} d\mu(z) \right) \left( \frac{1}{\text{diam}R_i} \int_{2B_i} N_i(z) d\mu(z) \right)^2 \\ &\quad\le C \frac{1}{\text{diam}R_i} \int_{2B_i} d\left(z, D_{p,t}\right) N_i(z)^{-2} d\mu(z).\end{aligned}$$

Therefore

$$\begin{aligned}&\sum_{i \in I(p,t)} \frac{\text{diam}R_i}{t} \times \frac{1}{t} \left( \frac{1}{\text{diam}R_i} \int_{2B_i} d\left(z, D_{p,t}\right)^{\frac{1}{3}} d\mu(z) \right)^3 \\ &\qquad\le \frac{C}{t^2} \sum_{i \in I(p,t)} \int_{2B_i} d\left(z, D_{p,t}\right) N_i(z)^{-2} d\mu(z) \\ &\qquad\le \frac{C}{t^2} \int_{\bigcup_{i \in I(p,t)} 2B_i} d\left(z, D_{p,t}\right) d\mu(z).\end{aligned}$$

We estimate this last quantity. If $i \in I(p,t)$, $R_i \cap B(p,t) \neq \emptyset$. If now $u \in B(p,t)$, $D(u) \le D(p) + t \le Ct$ so that there exists $u \in R_i$ such that $D(u) \le Ct$; hence $\text{diam}R_i \le Ct$, which implies $d\left(\pi(B_i), \pi(X)\right) \le Ct$. (We recall that $X$ is any point in $B(\tilde{X},t)\cap F$.) Hence, by Lemma 3.9 applied to a point $x$ in $2B_i$ (which satisfies $d(x)\leq 3\mathrm{diam} B_i\leq Ct$ ) and to the point $X$ which satisfies $d(X)\leq Ct$, we get $2B_i\subset B(X,Ct)$ so that, provided $k$ is large enough,

$$\begin{aligned}\frac{C}{t^2}\int_{\bigcup_{i\in I(p,t)}2B_i} d\left(z,D_{p,t}\right)d\mu(z) &\leq \frac{C}{t^2}\int_{B(X,Ct)} d\left(z,D_{p,t}\right)d\mu(z)\\ &\leq C\beta_1(X,t).\end{aligned}$$

Now, from estimates (4.1), (4.2), (4.3) and Lemmas 4.3 and 4.4 and because of the facts that $X$ is any point in $B(\tilde{X},t)\cap F$ and that $\mu(B(\tilde{X},t)\cap F)\geq \delta t$,

$$\begin{aligned}\gamma(p,t)^2 &\leq C\frac{1}{t}\int_{B(\tilde{X},t)\cap F}\beta_1(X,t)^2 d\mu(X)\\ &\quad +C\left\{\varepsilon\sum_{i\in I(p,t)}\left(\frac{\mathrm{diam} R_i}{t}\right)^2\right\}^2.\end{aligned}$$

We have then

$$\begin{aligned}\int_{U_0}\int_{\frac{D(p)}{60}}^2\gamma(p,t)^2\frac{dpdt}{t} &\leq C\int_{U_0}\int_{\frac{D(p)}{60}}^2\frac{1}{t}\int_{B(\tilde{X}(p,t),t)\cap F}\beta_1(X,t)^2 d\mu(X)\frac{dtdp}{t}\\ &\quad +C\varepsilon^2\int_{U_0}\int_{\frac{D(p)}{60}}^2\left\{\sum_{i\in I(p,t)}\left(\frac{\mathrm{diam} R_i}{t}\right)^2\right\}^2\frac{dtdp}{t}\\ &\leq C(a+b).\end{aligned}$$

We first look at the integral $a$. For any triple $(X,p,t)$ appearing in the computation, $|\pi(X)-p|\leq Ct$ and $\tilde{\delta}(X,t)\geq\frac{\delta}{C}$ (recall that the function $\tilde{\delta}$ appears in Definition 2.1) because Lemma 2.3 guarantees the existence of balls containing enough mass of $F$ in the ball $B(X,Ct)$. Hence

$$\begin{aligned}a &\leq \int_F\int_0^2 1\!\!1_{\{\tilde{\delta}(X,t)\geq\frac{\delta}{C}\}}\frac{1}{t}\left(\int_{p\in B(\pi(X),Ct)}dp\right)\beta_1(X,t)^2 d\mu(X)\frac{dt}{t}\\ &\leq C\int_F\int_0^2 1\!\!1_{\{\tilde{\delta}(X,t)\geq\frac{\delta}{C}\}}\beta_1(X,t)^2 d\mu(X)\frac{dt}{t}\\ &\leq Cc^2(\mu)\text{ by Corollary 2.4.}\end{aligned}$$

To estimate $b$, we remark that if $i\in I(p,t)$ then $\mathrm{diam} R_i\leq Ct$ (because $D(u)\leq D(p)+t\leq Ct$ on $B(p,t)$) so that

$$\sum_{i\in I(p,t)}\left(\frac{\mathrm{diam} R_i}{t}\right)^2\leq Ct\sum_{i\in I(p,t)}\frac{1}{t^2}\mathrm{diam} R_i\leq C.$$

Hence, noticing that $d(p, R_i) \leq t$ when $i \in I(p,t)$ and that $t \geq \frac{\mathrm{diam} R_i}{C}$ when $t \geq \frac{D(p)}{60}$, we obtain

$$\begin{aligned} b &\leq C\varepsilon^2 \int_{U_0} \int_{\frac{D(p)}{60}}^2 \sum_{i \in I(p,t)} \left(\frac{\mathrm{diam} R_i}{t}\right)^2 \frac{dtdp}{t} \\ &\leq C\varepsilon^2 \sum_{i \in I_0} (\mathrm{diam} R_i)^2 \int_{\frac{\mathrm{diam} R_i}{C}}^2 \int_{d(p,R_i) \leq t} \frac{dtdp}{t^3}. \end{aligned}$$

Therefore

$$\begin{aligned} b &\leq C\varepsilon^2 \sum_{i \in I_0} (\mathrm{diam} R_i)^2 \int_{\frac{\mathrm{diam} R_i}{C}}^2 (\mathrm{diam} R_i + t) \frac{dt}{t^3} \\ &\leq C\varepsilon^2 \sum_{i \in I_0} \mathrm{diam} R_i \\ &\leq C\varepsilon^2. \end{aligned}$$

Hence using these estimates, we get

$$\int_0^2 \int_{U_0} \gamma(p,t)^2 \frac{dpdt}{t} \leq C\varepsilon^2 + Cc^2(\mu) \leq C\varepsilon^2,$$

because $\eta \ll \varepsilon^2$, which ends the proof of Proposition 4.1. □

## 5. Calderón's formula and the size of $F_3$

From now on, we will extend $A$ to the whole line $D_0$ in a $C\alpha$-Lipschitz function of compact support. Let $\vec{D_0}$ be the line parallel to $D_0$ going through 0. Let $\nu : \vec{D_0} \to \mathbb{R}$ be an even, nonzero, $C^\infty$ function supported in $B(0,1)$ such that $\int_{\vec{D_0}} P\nu = 0$ for any affine function $P$ on $\vec{D_0}$.

We set $\nu_t(p_0) = \frac{1}{t}\nu(\frac{p_0}{t})$. When $f$ is a function defined on $D_0$, we write

$$(\nu_t * f)(p) = \int_{D_0} \nu_t(p-q) f(q) dq.$$

Calderón's formula (see [7, p.16,(5.9),(5.10)]) gives that, up to a normalization of $\nu$,

$$A(p) = \int_0^\infty (\nu_t * \nu_t * A)(p) \frac{dt}{t}.$$

We set $A = A_1 + A_2$ with

$$\begin{aligned} A_1 &= \int_2^\infty (\nu_t * \nu_t * A)(p) \frac{dt}{t} \\ &\quad + \int_0^2 (\nu_t * (\mathbb{1}_{D_0 \setminus U_0} (\nu_t * A)))(p) \frac{dt}{t} \\ A_2 &= \int_0^2 (\nu_t * (\mathbb{1}_{U_0} (\nu_t * A)))(p) \frac{dt}{t} \end{aligned}$$

where $U_0 = D_0 \cap B(0,10)$ as is defined just after Lemma 3.11.

LEMMA 5.1.
$$\int_D |\partial A_2|^2(p)dp \le C\int_0^2\int_{U_0}\gamma(p,t)^2\frac{dpdt}{t}.$$

*Proof.*
$$\partial A_2(p)=\int_0^2((\partial\nu)_t*(\mathbb{1}_{U_0}(\nu_t*A)))(p)\frac{dt}{t^2}.$$
We prove the $L^2$ estimate by a duality argument. For $F\in L^2(D)$,
$$\begin{aligned}\int_{D_0}F\partial A_2 &= \int_0^2\int_{D_0}\int_{D_0}F(p)(\partial\nu)_t(p-q)(\mathbb{1}_{U_0}(\nu_t*A))(q)\frac{dpdqdt}{t^2}\\ &= \int_0^2\int_{D_0}(tF*(\partial\nu)_t)(q)(\mathbb{1}_{U_0}(\nu_t*A))(q)\frac{dqdt}{t^3}\\ &\le \left(\int_0^2\int_{D_0}\mathbb{1}_{U_0}|\nu_t*A)|^2(q)\frac{dqdt}{t^3}\right)^{\frac12}\\ &\quad\times\left(\int_0^2\int_{D_0}|F*(\partial\nu)_t|^2(q)\frac{dqdt}{t}\right)^{\frac12}.\end{aligned}$$
By Plancherel's formula, using the fact that $\nu$ is radial, we have
$$\begin{aligned}\int_0^2\int_{D_0}|F*(\partial\nu)_t|^2(q)\frac{dqdt}{t} &= \int_0^2\int_{D_0}|\widehat{F}|^2(\xi)|(\widehat{\partial\nu})_t|^2(\xi)\frac{d\xi dt}{t}\\ &\le \int_{D_0}|\widehat{F}|^2(\xi)\int_0^\infty|(\widehat{\partial\nu})_t|^2(\xi)\frac{dt}{t}d\xi\\ &\le C(\nu)\int_{D_0}|F|^2(p)dp,\end{aligned}$$
so that
$$\int_{D_0}|\partial A_2|^2(p)dp\le C(\nu)\int_0^2\int_{D_0}\mathbb{1}_{U_0}|\nu_t*A|^2(q)\frac{dqdt}{t^3}.$$
Moreover, as the first moments of $\nu$ are zero, if $a$ is an affine function,
$$\begin{aligned}\left|\frac{\nu_t*A}{t}\right|(p) &= \left|\frac{\nu_t*(A-a)}{t}\right|(p)\\ &\le \frac1t\int_{B(p,t)\cap D_0}\left|\frac{\nu(\frac{p-q}{t})(A(q)-a(q))}{t}\right|dq\\ &\le \frac{C(\nu)}{t}\int_{B(p,t)\cap D_0}\left|\frac{A(q)-a(q)}{t}\right|dq,\end{aligned}$$
so that, taking the infimum over all affine functions $a$, we get
$$\left|\frac{\nu_t*A}{t}\right|(p)\le C(\nu)\gamma(p,t).$$

Hence we have

$\int_{D_0}|\partial A_2|^2(p)dp\le C(\nu)\int_0^2\int_{U_0}\gamma(p,t)^2\frac{dpdt}{t}.$ □

We set

$$U_1 = B(0,7) \cap D_0,$$
$$U_2 = B(0,4) \cap D_0.$$

LEMMA 5.2. *On $U_1$,*

$$|\partial A_1| \leq C\alpha,$$
$$|\partial^2 A_1| \leq C\alpha.$$

*Proof.* We set

$$\begin{aligned} A_1 &= A_{11} + A_{12} \text{ with} \\ A_{11}(p) &= \int_2^\infty (\nu_t * \nu_t * A)(p) \frac{dt}{t} \text{ and} \\ A_{12}(p) &= \int_0^2 (\nu_t * (\mathbb{1}_{D\setminus U_0}(\nu_t * A)))(p) \frac{dt}{t}. \end{aligned}$$

Note that on $U_1$, $A_{12}$ is zero for support reasons.

It remains to estimate $A_{11}$. We set

$$\psi = \int_2^\infty \nu_t * \nu_t \frac{dt}{t}.$$

Then

$$\begin{aligned} A_{11} &= \psi * A, \\ \partial A_{11} &= \psi * \partial A, \\ \partial^2 A_{11} &= \partial\psi * \partial A, \end{aligned}$$

so that

$$\begin{aligned} \|\partial A_{11}\|_\infty &\leq \|\partial A\|_\infty \int |\psi| \text{ and} \\ \|\partial^2 A_{11}\|_\infty &\leq \|\partial A\|_\infty \int |\partial\psi|. \end{aligned}$$

As it is known that $\|\partial A\|_\infty \leq C\alpha$, we only have to evaluate $\int |\psi|$ and $\int |\partial\psi|$.

$$\begin{aligned} \int_{|p|\leq 10} |\psi(p)| dp &\leq \int_2^\infty \int_{|p|\leq 10} \frac{1}{t} |\nu|(\frac{p}{t}) \int |\nu|(\frac{q-p}{t}) \frac{dq}{t} dp \frac{dt}{t} \\ &\leq \left(\int |\nu|\right) \|\nu\|_\infty \int_2^\infty \int_{|p|\leq 10} dp \frac{dt}{t^2} \\ &\leq C(\nu). \end{aligned}$$

Moreover,

$$\psi = \delta_0 - \int_0^2 \nu_t * \nu_t \frac{dt}{t}$$

because of Calderón's formula so that, as $\nu$ is zero outside $B(0,1)$, $\psi(p)=0$ if $|p|\geq 10$ ($\operatorname{supp}(\nu_t*\nu_t)\subset B(0,2t)\subset B(0,4)$ for $0\leq t\leq 2$).

We can do the same for $\int|\partial\psi|$ and this ends the proof of Lemma 5.2. □

Define the maximal function

$$N(A_2)(p)=\sup_B\left\{\frac{1}{|B|}\int_B\frac{|A_2-m_BA_2|}{|B|}\right\},$$

where the supremum is over all balls $B$ containing $p$ of radius $\leq 2$. Now we may state:

LEMMA 5.3.

$$\int_{D_0}N(A_2)^2\leq C\int_{D_0}|\partial A_2|^2.$$

*Proof.* By Poincaré's inequality,

$$\frac{m_B(|A_2-m_BA_2|)}{r}\leq Cm_B(|\partial A_2|)$$

so that

$$N(A_2)(p)\leq C\sup_{p\in B}m_B(|\partial A_2|);$$

hence, by the Hardy-Littlewood maximal inequality,

$$\int_{D_0}N(A_2)^2(p)dp\leq C\int_{D_0}|\partial A_2|^2(p)dp.$$ □

LEMMA 5.4. *Set* $\operatorname{osc}_B A_2=\sup_{p\in B}|A_2(p)-m_BA_2|$ *and let* $r$ *be the radius of* $B$. *Then, if* $B\subset U_1$,

$$\operatorname*{osc}_B A_2\leq Cr\left\{\frac{m_B(|A_2-m_BA_2|)}{r}\right\}^{\frac{1}{2}}\alpha^{\frac{1}{2}}.$$

*Proof.* Let $B\subset U_1$ and set $\lambda=\operatorname{osc}_B A_2=|A_2(q)-m_BA_2|$ for a point $q\in B$. As $\|\partial A_2\|_{L^\infty(B)}\leq C\alpha$, $|A_2(p)-m_BA_2|\geq\frac{\lambda}{2}$ when $p\in B$ and $d(p,q)\leq\frac{\lambda}{2C\alpha}$.

- If $\frac{\lambda}{2C\alpha}\leq r$,

$$\int_B|A_2(p)-m_BA_2|dp\geq\frac{\lambda}{2C}\frac{\lambda}{2C\alpha};$$

hence

$$\lambda^2\leq Cr^2\alpha\frac{m_B(|A_2(p)-m_BA_2|)}{r}.$$

- If $\frac{\lambda}{2C\alpha} \geq r$, $|A_2(p) - m_B A_2| \geq \frac{\lambda}{2}$ on more than half of the ball $B$ so that

$$m_B(|A_2 - m_B A_2|) \geq \frac{\lambda}{C}.$$

Moreover, by Poincaré's inequality,

$$\begin{aligned} \frac{m_B(|A_2 - m_B A_2|)}{r} &\leq C m_B(|\partial A_2|) \\ &\leq C\|\partial A_2\|_{L^\infty(B)} \\ &\leq C\alpha, \end{aligned}$$

so that, summarizing these inequalities, we get the result. □

Lemma 5.5. *For a number $\theta > 0$, set $H_\theta = \{p \in U_2, N(A_2)(p) \leq \theta^2\alpha\}$. If $B = B(p_0, r)$ intersects $H_\theta$ and $r \leq \theta$, then*

$$\sup_{p\in B} |A(p) - \{A(p_0) + \partial A_1(p_0)(p - p_0)\}| \leq Cr\theta\alpha.$$

*Proof.* If $p \in B$,

$$\begin{aligned} &|A(p) - \{A(p_0) + \partial A_1(p_0)(p - p_0)\}| \\ &\quad \leq |A_2(p) - A_2(p_0))| + |A_1(p) - \{A_1(p_0) + \partial A_1(p_0)(p - p_0)\}| \\ &\quad \leq 2 \operatorname*{osc}_B A_2 + C\alpha r^2 \text{ (by Taylor and Lemma 5.2)} \\ &\quad \leq Cr \left\{ \frac{m_B(|A_2 - m_B A_2|)}{r} \right\}^{\frac{1}{2}} \alpha^{\frac{1}{2}} + C\alpha r^2 \\ &\quad \leq Cr(N(A_2)(p_1))^{\frac{1}{2}} \alpha^{\frac{1}{2}} + C\alpha r^2 \text{ (where } p_1 \in B) \\ &\quad \leq Cr\theta\alpha \text{ by taking } p_1 \in H_\theta \cap B. \end{aligned}$$

□

If $\Delta_B$ is the line which is the graph of the function $p \mapsto A(p_0) + \partial A_1(p_0)(p - p_0)$,

$$\sup_{x\in\Gamma\cap\pi^{-1}(B)} \frac{d\,(x, \Delta_B)}{r} \leq C\theta\alpha.$$

Lemma 5.6. *If $\theta > 0$ is given, there exists $\varepsilon_0 > 0$ such that if $\varepsilon < \varepsilon_0$ then* $\operatorname{angle}(D_{x,t}, D_0) \leq \frac{\alpha}{100}$ *for any $(x, t) \in S$, $100t \geq \theta$.*

*Proof.* Let $(x, t)$ be such a couple and let $k$ be such that $2^{(k+1)} \leq t \leq 2^k$. We have, by Lemma 2.6,

$$\operatorname{angle}(D_{x,t}, D_{x,2t}) \leq C\varepsilon,$$

so that

$$\operatorname{angle}(D_{x,t}, D) \leq \sum_{j=0}^{k} \operatorname{angle}(D_{x,2^j t}, D_{x,2^{j+1}t}),$$

so that

$$\text{angle}(D_{x,t}, D) \leq (k+1)C\varepsilon \leq \frac{\alpha}{100},$$

provided $t \geq \frac{\theta}{100}$ and $\varepsilon$ is chosen after $\theta$. □

We set

$$\tilde{\tilde{F}} = \left\{ x \in \tilde{F}, \forall t \in (0,2), \mu(\tilde{F} \cap B(x,t)) \geq \frac{99}{100}\mu(F \cap B(x,t)) \right\}.$$

LEMMA 5.7. $\mu(F \backslash \tilde{\tilde{F}}) \leq C\varepsilon^{\frac{1}{2}}$.

*Proof.* It is enough to evaluate $\mu(\tilde{F} \backslash \tilde{\tilde{F}})$ because we already know how to evaluate $\mu(F \backslash \tilde{F})$. Now

$$\tilde{F} \backslash \tilde{\tilde{F}} \subset \bigcup_{x \in \tilde{F} \backslash \tilde{\tilde{F}}} B(x, t_x) \cap \tilde{F},$$

where $B(x,t_x)$ satisfies $\mu(B(x,t_x) \cap \tilde{F}) \leq 100\mu(B(x,t_x) \cap (F \backslash \tilde{F}))$. Hence, by Besicovitch's covering lemma, we get families $\mathcal{B}_n$, $n = 1, \ldots, N$, of disjoint balls $B(x,t_x)$, whose union is still a covering of $\tilde{F} \backslash \tilde{\tilde{F}}$ so that

$$\begin{aligned}
\mu(\tilde{F} \backslash \tilde{\tilde{F}}) &\leq \sum_{n=1}^{N} \sum_{B \in \mathcal{B}_n} \mu(B \cap \tilde{F}) \\
&\leq 100 \sum_{n=1}^{N} \sum_{B \in \mathcal{B}_n} \mu(B \cap (F \backslash \tilde{F})) \\
&\leq 100N\mu(F \backslash \tilde{F}) \\
&\leq C\varepsilon^{\frac{1}{2}}.
\end{aligned}$$

□

LEMMA 5.8. *If* $x \in F_3 \cap \tilde{\tilde{F}}$, $d(\pi(x), H_\theta) > h(x)$ *and* $h(x) \leq \frac{\theta}{100}$.

*Proof.* Let us recall first that $d(x) \leq h(x)$ (because of Remark 3.8) and that $h(x) \leq \frac{\theta}{100}$ because of Remark 3.3. Suppose that $d(\pi(x), H_\theta) < h(x)$. Setting $B = B(x, 2h(x))$, we would have $\pi(B) \cap H_\theta \neq \emptyset$, if $x' \in B \cap \tilde{F}$ (which is the case for 99 percent of $x$ in $B \cap F$),

$$\begin{aligned}
d(x', \pi(x') + A(\pi(x'))) &\leq \varepsilon^{\frac{1}{2}} d(x') \\
&\leq \varepsilon^{\frac{1}{2}}(d(x) + 2h(x)) \\
&\leq 3\varepsilon^{\frac{1}{2}} h(x),
\end{aligned}$$

so that

$$\begin{aligned} d\,(x',\Delta_B) &\leq d\,(x',\pi(x')+A(\pi(x'))) + d\,(\pi(x')+A(\pi(x')),\Delta_B) \\ &\leq (3\varepsilon^{\frac{1}{2}}+C\theta\alpha)h(x) \\ &\leq C\theta\alpha h(x). \end{aligned}$$

Hence $\text{angle}(D_{x,h(x)},\Delta_B) \leq \frac{\alpha}{100}$.

We may apply the same argument with the ball $B' = B(x,\frac{\theta}{100})$ and get that

$$\text{angle}(D_{x,\frac{\theta}{100}},\Delta_{B'}) \leq \frac{\alpha}{100}.$$

But $\Delta_B = \Delta_{B'}$ because these lines only depend on the projection of the center of the ball; hence

$$\text{angle}(D_{x,h(x)},D_{x,\frac{\theta}{100}}) \leq \frac{\alpha}{50}.$$

Moreover, by Lemma 5.6,

$$\text{angle}(D_0,D_{x,\frac{\theta}{100}}) \leq \frac{\alpha}{50}$$

so that

$$\text{angle}(D_0,D_{x,h(x)}) \leq \frac{\alpha}{25},$$

which is impossible because of Remark 3.3. □

PROPOSITION 5.9. *Provided the parameters $\theta$, $\alpha$ and $\varepsilon$ are well chosen,*

$$\mu(F_3) \leq 10^{-6}.$$

*Proof.* We only have to evaluate $\mu(F_3\cap\tilde{\tilde{F}})$, as follows:

$$F_3\cap\tilde{\tilde{F}} \subset \bigcup_{F_3\cap\tilde{\tilde{F}}} B(x,2h(x))\cap\tilde{\tilde{F}}.$$

Now by Besicovitch's covering lemma and the upper control on $\mu$, we get

$$\begin{aligned} \mu(F_3\cap\tilde{\tilde{F}}) &\leq \sum_{n=0}^{N}\sum_{B\in\mathcal{B}_n}\mu(B\cap\tilde{\tilde{F}}) \\ &\leq C_0\sum_{n=0}^{N}\sum_{B\in\mathcal{B}_n}\text{diam}B, \end{aligned}$$

where the balls $B$ are of the type $B(x,2h(x))$ for a point $x\in F_3\cap\tilde{\tilde{F}}$ and where two balls of the same family $\mathcal{B}_n$ are disjoint.

If $B$ and $B'$ are two balls of the same family $\mathcal{B}_n$, provided $\varepsilon^{\frac{1}{2}}$ is very small compared to $\alpha$, the line going through the centers of $B$ and $B'$ has slope $\leq C\alpha$. This is because the center of $B$ is at distance less than $\varepsilon^{\frac{1}{2}}\text{diam}B$ from the graph of $A$ (see Fig. 2, above) which is a $C\alpha$-Lipschitz function and the same is true for $B'$. Hence, provided $\alpha$ is very small, the projections of $\frac{1}{2}B$ and $\frac{1}{2}B'$ on $D$ are disjoint. We know that the projections of these balls do not meet $H_\theta$ and are contained in $U_2$, so that

$$\sum_{B \in \mathcal{B}_n} \mathrm{diam} B \leq 2\mu(U_2 \backslash H_\theta),$$

which implies

$$\mu(F_3 \cap \tilde{\tilde{F}}) \leq 2C_0 N \mu(U_2 \backslash H_\theta).$$

Now, by Lemma 5.3, Lemma 5.1 and Proposition 4.1,

$$\int_{D_0} N(A_2)^2 \leq C\varepsilon^2,$$

so that, from the definition of $H_\theta$,

$$\mu(U_2 \backslash H_\theta) \leq C \frac{\varepsilon^2}{\theta^4 \alpha^2}.$$

Hence, choosing $\varepsilon$ after $\theta$ and $\alpha$, we will get

$$\begin{aligned} \mu(F_3 \cap \tilde{\tilde{F}}) &\leq 10^{-7} \text{ and} \\ \mu(F \backslash \tilde{\tilde{F}}) &\leq 10^{-7}, \end{aligned}$$

which gives the proposition. □

Université Paris XI, Orsay, France
*E-mail address*: Jean-Christophe.LEGER@math.u-psud.fr

References

[1] G. David, Unrectifiable sets have vanishing analytic capacity, Revista Mat. Iberoamericana **14** (1998), 369–479.

[2] G. David and P. Mattila, Removable sets for Lipschitz harmonic functions in the plane, Preprint, Université Paris-Sud, 1997.

[3] G. David and S. Semmes, *Singular Integrals and Rectifiable Sets in* $\mathbb{R}^n$: *Au delà des Graphes Lipschitziens*, No. 193 in Astérisque, SMF, 1991.

[4] K. J. Falconer, *Geometry of Fractals Sets*, Cambridge University Press, 1984.

[5] P. W. Jones, Rectifiable sets and the traveling salesman problem, Invent. Math. **102** (1990), 1–15.

[6] M. S. Melnikov, P. Mattila, and J. Verdera, The Cauchy integral, analytic capacity, and uniform rectifiablity, Ann. of Math. **144** (1996), 127–136.

[7] Y. Meyer, *Ondelettes et Opérateurs* I: *Ondelettes* (French), Actualites Math. Paris: Herman, Editeurs des Sciences et des Arts, 1990.

[8] A. G. Vithushkin, The analytic capacity of sets in problems of approximation theory, Uspekhi Mat. Nauk **22** (1967), 141–199; English translation in Russian Math. Surveys **22** (1967), 139–200.